\documentclass{article}

\usepackage[round]{natbib}
\usepackage{amsmath, amssymb, amsfonts}
\usepackage{amsthm, thmtools, thm-restate}
\usepackage{mathtools}
\usepackage{enumitem}

\usepackage{algpseudocode}

\usepackage{epsfig, graphicx, subfigure}
\usepackage{bbold}

\usepackage{upgreek}

\usepackage{tikz}
\usetikzlibrary{positioning}
\usetikzlibrary{shapes.geometric, shapes}
\usetikzlibrary{fit}
\usepackage{stackengine}
\def\shft#1{\stackunder[43pt]{}{\kern112pt #1}}

\tikzset{
	module/.style={draw,thick,  shape=rectangle, rounded corners = 0.5ex, minimum width =2em, minimum height =2em},
	bigmodule/.style={draw, thick, shape=rectangle, rounded corners = 1ex, align=right, minimum width =3em, minimum height =3em, inner sep=2ex},
	square/.style={regular polygon,regular polygon sides=4}
}

\usepackage{color}
\definecolor{puorange}{rgb}{0.80,0.20,0}
\definecolor{bluegray}{RGB}{65,105,225}
\definecolor{greengray}{rgb}{0.05,0.50,0.15}
\definecolor{darkbrown}{rgb}{0.40,0.2,0.05}
\definecolor{darkcyan}{RGB}{0,160,200}
\definecolor{black}{rgb}{0,0,0}
\definecolor{darkgray}{rgb}{0.3,0.3,0.3}
\definecolor{darkmagenta}{RGB}{139,0,139}
\definecolor{darkpurple}{RGB}{128,0,128}
\definecolor{darkred}{RGB}{90,0,0}
\definecolor{nicered}{RGB}{178,34,34}
\definecolor{darkorange}{RGB}{255,90,0}
\definecolor{darkblue}{RGB}{20,20,100}
\definecolor{verydarkorange}{RGB}{180,50,0}
\definecolor{teal}{RGB}{0,128,128}
\definecolor{indigo}{RGB}{75,0,130}
\usepackage[colorlinks,citecolor=bluegray,linkcolor=darkred,urlcolor=darkcyan,breaklinks]{hyperref}

\newtheorem{theorem}{Theorem}[section]

\newtheorem{lemma}[theorem]{Lemma}
\newtheorem{corollary}[theorem]{Corollary}

\newtheorem{fact}[theorem]{Fact}

\newcommand{\reals}{{\mathbb R}}

\newcommand{\ones}{\operatorname{\mathbf 1}}
\newcommand{\idm}{\operatorname{Id}}

\DeclareMathOperator*{\argmax}{argmax}
\DeclareMathOperator*{\argmin}{argmin}
\newcommand{\bigO}{O}
\newcommand{\dom}{\operatorname{dom}}

\def\shortdisplay{\setlength{\abovedisplayskip}{5pt}%
	\setlength{\belowdisplayskip}{5pt}%
	\setlength{\abovedisplayshortskip}{2pt}%
	\setlength{\belowdisplayshortskip}{2pt}}

\let\oldselectfont\selectfont
\def\selectfont{\oldselectfont\shortdisplay}

\newcommand{\point}{x_0}
\newcommand{\labl}{y}

\newcommand{\obj}{h}
\newcommand{\reg}{\omega}
\newcommand{\loss}{\mathcal{L}}

\newcommand{\ctrl}{w}
\newcommand{\ctrls}{{\ctrl}}
\newcommand{\state}{x}
\newcommand{\states}{\state}

\newcommand{\costate}{\lambda}
\newcommand{\costates}{\costate}

\newcommand{\chain}{f}
\newcommand{\chainnet}{\psi}

\newcommand{\dyn}{\phi}

\newcommand{\horizon}{\tau}
\newcommand{\layer}{t}

\newcommand{\auxctrl}{v}

\newcommand{\auxstate}{y}

\newcommand{\dimctrl}{p}
\newcommand{\dimctrls}{{\dimctrl}}

\newcommand{\dimin}{d}
\newcommand{\dimout}{m}

\newcommand{\dimvar}{d}
\newcommand{\diminter}{k}

\newcommand{\lip}{\ell}
\newcommand{\smooth}{L}

\newcommand{\stepsize}{\gamma}
\newcommand{\lin}{\ell}

\newcommand{\lag}{\mathcal{L}}

\newcommand{\prox}[1]{\operatorname{prox}(#1)}

\newcommand{\target}{z}

\newcommand{\reglag}{\kappa}

\newcommand{\diffprog}{\mathcal{P}}
\newcommand{\diffprogfunc}{f}
\newcommand{\morprog}{\mathcal{M}}
\newcommand{\env}{\operatorname{env}}

\newcommand{\BP}{\operatorname{BP}}
\newcommand{\GBP}{\operatorname{GBP}}
\newcommand{\MBP}{\operatorname{MBP}}
\newcommand{\IBP}{\operatorname{IBP}}

\newcommand{\GD}{\operatorname{GD}}

\newcommand{\infunc}{f}
\newcommand{\outfunc}{h}
\newcommand{\anv}{\operatorname{anv}}

\newcommand{\adjvar}{\lambda}
\newcommand{\dir}{\mu}
\newcommand{\relu}{\operatorname{relu}}
\newcommand{\aug}{\operatorname{aug}}
\newcommand{\pen}{\mathcal{P}}

\newcommand{\invsign}{\dagger}
\newcommand{\augparam}{\kappa}
\newcommand{\pgstepsize}{\rho}
\newcommand{\scaling}{\sigma}

\newcommand{\auxstepsize}{\alpha}
\newcommand{\auxxstepsize}{\beta}
\usepackage{algorithm, algorithmicx}
\usepackage{natbib}
\usepackage{adjustbox}
\usepackage{fullpage}

\let\originalparagraph\paragraph
\renewcommand{\paragraph}[2][.]{\originalparagraph{#2#1}}

\usepackage{caption}
\title{Differentiable Programming \` a la Moreau}
\date{}
\author{Vincent Roulet, Zaid Harchaoui\\ 
	Department of Statistics, University of Washington, Seattle, USA}

\begin{document}
\maketitle
\begin{abstract}
	The notion of a Moreau envelope is central to the analysis of first-order optimization algorithms for machine learning. Yet, it has not been developed and extended to be applied to a deep network and, more broadly, to a machine learning system with a differentiable programming implementation. We define a compositional calculus adapted to Moreau envelopes and show how to integrate it within differentiable programming. The proposed framework casts in a mathematical optimization framework several variants of gradient back-propagation related to the idea of the propagation of virtual targets
\end{abstract}

\section*{Introduction}
We consider objectives  characterized by a function $\chain$ decomposable in a sequence of elementary operations $\dyn_\layer$, which,
given an initial point $\point$,  maps a sequence of variables $\ctrls = (\ctrl_1, \ldots,\ctrl_\horizon)$ to an output
\begin{align}\label{eq:chain}
	\chain(\ctrls)  &= \state_\horizon,  \\
	\mbox{s.t.} \ \state_\layer & =	\dyn_\layer(\ctrl_\layer, \state_{\layer{-}1}) \ \mbox{for} \ \layer = 1, \ldots, \horizon.\nonumber
\end{align}
Such a dynamical structure typically arises in deep learning problems, where $\dyn_\layer$ are layers and $\ctrl_1, \ldots,\ctrl_\horizon$ are the weights of all layers, and in nonlinear discrete control problems, where $\dyn_\layer$ are nonlinear dynamics and $\ctrl_1, \ldots,\ctrl_\horizon$ represent a sequence of controls~\citep{zhang2021dive, slotine1991applied}. Given a dynamical structure~\eqref{eq:chain}, the optimization problem then consists in solving 
$
\min_{\ctrls} \obj(\chain(\ctrls)) 
$
for $\obj$ a cost on the output of the dynamical system. 

Standard gradient-based optimization methods can be used to solve such problems. Obtaining the gradient then amounts to applying the chain-rule, which is nowadays usually implemented using automatic differentiation for deep networks and other complex models in a differentiable programming framework~\citep{paszke2017automatic, tensorflow2015-whitepaper, bolte2020mathematical}. 
As differentiable programming stands out as a computational framework tailored for training models using first-order optimization, one may ask how the notion of Moreau envelope could fit into it and expand its scope. 

Indeed, the notion of Moreau envelope~\citep{moreau1962fonctions, yosida2012functional,  martinet1970regularisation, martinet1972determination,rockafellar1976monotone, attouch1977convergence}  has arisen as a central notion in the analysis of first-order optimization algorithms for machine learning~\citep{duchi:ruan:2018, lin2018catalyst,drusvyatskiy2019efficiency}. To blend Moreau envelopes into differentiable programming, one needs to define a calculus adapted to Moreau envelopes. We propose a framework to define such a calculus and show how to integrate it within differentiable programming. We show how previous proposals of smoother alternatives to gradient back-propagation fit into our framework. 
We present numerical results in deep learning and nonlinear control.

\paragraph{Related work}
The computational building blocks we consider are similar to the ones considered in variants of gradient back-propagation, which can be traced back to the now called target propagation~\citep{lecun1988theoretical, lecun1989gemini, rohwer1990moving,mirowski2009dynamic}.  Target propagation can be described as using approximate inverses of layers when computing an update of a deep network~\citep{lee2015difference, meulemans2020theoretical,  ahmad2020gait}. The moving targets that minimize the overall objective are back-propagated via approximate layer inverses. The layer weights are then updated by minimizing the distance between the output of the layer and the given moving target. These algorithms were found to be effective in some settings and were, for the most part, motivated by empirical observations. Penalized formulations of the training problem have also been considered to decouple the optimization of the weights in a distributed way \citep{carreira2014distributed, taylor2016training, gotmare2018decoupling}. 
Finally, our framework encompasses the proximal back-propagation algorithm of~\citep{frerix2018proximal} which mixes the classical gradient back-propagation and a proximal step to update the weights of a deep forward network, to get a proximal-type gradient back-propagation.

\paragraph{Notations}
For a function $f: \reals^ \dimin \rightarrow\reals^ \dimout$, we denote $\nabla f(x) = (\partial_{x_i} f_j(x) )_{i\in \{1, \ldots, \dimin\}, j \in \{1, \ldots, \dimout\}}\in \reals^{\dimin \times \dimout}$ the gradient of $f$ at $x$, i.e., the transpose of the Jacobian.  We denote the Lipschitz continuity parameter of $f$ as
$\lip_f = \sup_{\substack{x, y \in \reals^ \dimin,  x\neq y}}  \|f(x) {-} f(y)\|_2 /\|x-y\|_2$ 
and the smoothness parameter of $f$, i.e., the Lipschitz continuity parameter of its gradient, as  $\smooth_f = \sup_{\substack{x, y \in \reals^ \dimin,  x\neq y}}  \|\nabla f(x) - \nabla f(y)\|_{2}/\|x-y\|_2$,
where $\|A\|_2$ denotes the spectral norm of a matrix $A$.

\section{Differentiable Program for the Moreau Envelope}\label{sec:moreau}
Key to the minimization of dynamical systems of the form~\eqref{eq:chain}   is the availability of first-order information via automatic differentiation in a differentiable programming framework.
Formally, a differentiable program $\diffprog$ implements the evaluation of a function $\diffprogfunc$ and enables the computation of any gradient-vector product on the evaluated point. Formally, $\mathcal{P}$ can be defined as
	\[
	\diffprog : \begin{cases}
	\reals^\dimin & \rightarrow \reals^\dimout \times (\reals^\dimout \rightarrow \reals^\dimin)\\
w & \mapsto (\diffprogfunc(w), \lambda \mapsto \nabla \diffprogfunc(w) \lambda),
	\end{cases}
	\]
such that in practice, in, e.g., PyTorch~\citep{paszke2017automatic}, $\mathcal{P}$ consists in evaluating $\texttt{out  = f(w)}$ such that gradient-vector products are accessed as $ \nabla \diffprogfunc(w) \lambda = \texttt{autograd(out, w, lambda)}$. 
To minimize an objective of the form $h(f(w))$ with a gradient descent, it suffices then to compute $f(w)$ through the program $\mathcal{P}$ and access $\nabla (h\circ f)(w) = \nabla f(w) \nabla h(w)$ to perform an update of the form $w \leftarrow w - \stepsize \nabla f(w) \nabla h(w)$ for some $\stepsize>0$. 

The computation of the gradient-vector product can be seen as the minimization of a regularized linear approximation of the objective, i.e., 
$
\nabla f(w) \lambda = \argmin_{v\in \reals^\dimin}\lambda^\top \lin_{f}^w(w-v) +\|v\|_2^2/2,
$
where $\lin_{f}^w(v) = f(w) + \nabla f(w)^\top (v-w)$ is the linear approximation of $f$ around $w$ that can be decomposed into the elementary computations $\dyn$ defined in~\eqref{eq:chain}. This viewpoint serves as a basis for efficient implementations of Newton or Gauss-Newton steps by decomposing the quadratic approximations of the objective into the elementary computations and solving the resulting subproblem by dynamic programming as done in nonlinear control~\citep{dunn1989efficient, wright1991partitioned} or for generic graphs of computations~\citep{srinivasan2021computing}. Following this approach, we seek to take advantage of the decomposition of a function $f$ into elementary computations to compute an oracle on the objective based on its Moreau envelope.

\subsection{Moreau Envelope}
\paragraph{Moreau gradient}
The Moreau envelope  defines  an oracle through the minimization of the function rather than using an approximation of the function~\citep{moreau1962fonctions, bauschke2011convex}. Formally, for a real function $f:\reals^\dimin \rightarrow \reals$, $w\in \reals^d$ and $\bar \stepsize>0$ such that\footnote{The Moreau envelope of $\stepsize f$ is guaranteed to exist for $0\leq \stepsize <\bar \stepsize$ as long as  $\bar \stepsize$ is such that $\inf_{v \in \reals^\dimin}  \bar \stepsize f(w-v) + \|v\|_2^2/2$ is finite. We consider condition~\eqref{eq:mor_cond} to simplify the exposition of the smoothness properties of the Moreau gradient. Note that convexity of $v \mapsto  \bar \stepsize f(w-v) + \|v\|_2^2/2$ on a sufficiently large ball around the origin is sufficient to guarantee that $\inf_{v \in \reals^\dimin}   \stepsize f(w-v) + \|v\|_2^2/2$ is finite for $\stepsize< \bar \stepsize$. }
\begin{equation}\label{eq:mor_cond}
v \mapsto  \bar \stepsize f(w-v) + \|v\|_2^2/2 \  \mbox{is convex},
\end{equation}
the Moreau envelope of $\stepsize f$ on $w \in \reals^d$ for $0 \leq \stepsize <\bar \stepsize$ is defined as  
\begin{align}\label{eq:moreau}
\env(\stepsize f)(w)  & = \inf_{v \in \reals^\dimin}  \left\{\stepsize f(w-v) + \|v\|_2^2/2\right\}.
\end{align}
The gradient of the Moreau envelope of $\stepsize f$ on $w \in \reals^d$, called hereafter the Moreau gradient of $f$ on $w\in \reals^d$ with parameter $\stepsize$, is defined as
\begin{align*}
 \nabla \env(\stepsize f)(w) 
&  =  \argmin_{v \in \reals^\dimin} \left\{\stepsize f(w-v) + \|v\|_2^2/2 \right\}.
\end{align*}
While the gradient of $f$ may not be Lipschitz-continuous or may not even be defined, the Moreau gradient is well-defined for any  $0 \leq \stepsize <\bar \stepsize$ and is $\max\{\stepsize/(\bar \stepsize- \stepsize), 1\}$-Lipschitz-continuous~\citep[Corollary 3.4]{hoheisel2020regularization}, that is,  $1$-Lipschitz continuous for any $0\leq \stepsize \leq \bar \stepsize/2$.
On the other hand, the Moreau gradient defines a first-order optimality condition for minimizing $f$, since if $w^*\in \argmin_{w \in \reals^d} f(w)$, then $\nabla \env(\stepsize f)(w^*) = 0$~\citep{bauschke2011convex}. 

Equipped with an oracle that gives access to the Moreau gradients, we can define a Moreau gradient descent with step-size $\stepsize$ by the updates  $w \leftarrow w - \nabla \env(\stepsize f)(w)$. 
The parameter $\stepsize$ acts as a stepsize for the oracle: the smaller the $\stepsize$, the smaller the Moreau gradient. Note that the stepsize is part of the definition of the Moreau gradient such that  $\stepsize \nabla \env(f)(w) \neq \nabla \env(\stepsize f)(w)$.
In the vocabulary of proximal operators, a Moreau gradient descent is equivalent to a proximal point method $w \leftarrow \prox{\stepsize f}(w)$, where the proximal operator is defined as $\prox{\stepsize f}(w) = \argmin_{v \in \reals^\dimin} f(v) + \|w-v\|_2^2/(2\stepsize) = w - \nabla \env(\stepsize f)(w)$. Convergence proofs of proximal point methods for any convex function have been shown by~\citet{rockafellar1976monotone, bauschke2011convex} by using that exact proximal point iterations  ensure decreasing objective values. Here we focus on Moreau gradients which define smooth surrogates of the gradient of $f$ at a cost of solving~\eqref{eq:moreau}.

\paragraph{Approximate Moreau gradient}
In practice, one usually approximates the Moreau envelope using an optimization algorithm; see,~{e.g.},~\cite{lin2018catalyst}. Namely, for  $f$ differentiable, $w \in \reals^d$ and $0 \leq \stepsize < \bar \stepsize$ such that condition~\eqref{eq:mor_cond} is satisfied, the Moreau gradient can be computed as
\begin{align}\nonumber
\nabla \env(\stepsize f)(w) = \lim\limits_{k \rightarrow +\infty}\mathcal{A}_k\left(\stepsize f(w-\cdot) + \|\cdot\|_2^2 /2\right),
\end{align}
where $\mathcal{A}_k\left(h \right)$ is the $k\textsuperscript{th}$ output of an algorithm $\mathcal{A}$, such as gradient descent,  applied to minimize a function $h$. For example, if condition~\eqref{eq:mor_cond} is satisfied and $f$ is smooth, a gradient descent can estimate the Moreau gradient up to $\varepsilon$ accuracy at a logarithmic cost in~$\varepsilon$~\citep{lin2018catalyst}. We are then interested in developing optimization algorithms that utilize approximate Moreau gradients such as an approximate Moreau gradient descent defined by the iterates
\begin{equation}\label{eq:mgd}
	w^{(k+1)} = w^{(k)} - \widehat \nabla \env (\stepsize f)(w^{(k)}),
\end{equation}
where $\widehat \nabla \env (\stepsize f)(w) \approx \nabla \env(\stepsize f)(w)$ is an approximation of the Moreau gradient of $f$ at $w$ with stepsize $\stepsize$. An overview of the convergence guarantees of such an approach is provided in Appendix~\ref{app:moreau}.

\paragraph{Augmented Moreau gradient}
Up to a change of variables and a rescaling, the classical Moreau envelope can be written as $\min_{v \in \reals^\dimin} f(v) + \|w-v\|_2^2/(2\stepsize)$, i.e., a minimization of $f$ with a regularization term forcing the minimum to be close to the current point. We additionally consider an augmented definition of the Moreau envelope that regularizes the minimum both in the input and output space, i.e., $\min_{v \in \reals^\dimin} f(v) + \|w-v\|_2^2/(2\stepsize_1) + (f(w) - f(v))^2/(2\stepsize_2)$. 
For ease of presentation, we consider first a simplified definition of an augmented Moreau envelope of $f$ on $w\in \reals^\dimin$ with parameter $\stepsize\geq 0$ 
\[
\anv(f)(w; \stepsize)  = \inf_{v\in \reals^\dimin} \left\{\stepsize f(w-v) + \frac{1}{2} \|v\|_2^2 + \frac{1}{2} (f(w-v)-f(w))^2\right\}.
\]
The augmented Moreau envelope is well defined for any $\stepsize\geq 0 $ as we have $\anv( f)(w; \stepsize) \geq \stepsize f(w) - \stepsize^2/2$ for any $w \in \reals^\dimin$. Moreover if $f$ is such that  $v\mapsto  (f(w-v)-f(w))^2/2 + \|v\|_2^2/4$ is convex, the minimizer is unique for  $\stepsize \leq \bar \stepsize/2$ with $\bar \stepsize$ satisfying~\eqref{eq:mor_cond} and defines an 
augmented Moreau gradient of $f$ on $w\in \reals^\dimin$ with parameter $\stepsize$ as 
\begin{align}
	\nabla \anv( f)(w; \stepsize) 
	& 	= \argmin_{v \in \reals^\dimin}\left\{\stepsize f(w-v) + \frac{1}{2} \|v\|_2^2 + \frac{1}{2} (f(w-v)-f(w))^2\right\}, \nonumber\\
	& =  \argmin_{v \in \reals^\dimin}\left\{ (f(w-v) - f(w) + \stepsize)^2 +  \|v\|_2^2 \right\}. \label{eq:aug_mor_grad}
\end{align}
Here, for ease of presentation, we use the symbol $\nabla \anv( f)(w; \stepsize)$ to denote the augmented Moreau gradient, though if the function $g:w\mapsto \anv(f)(w;\stepsize)$ is differentiable, its gradient is not $\nabla \anv( f)(w; \stepsize)$.
The augmented Moreau gradient can then be interpreted as a regularized inversion of $f$ to ensure a decrease of $\stepsize$ in objective values. More generally,  we define an augmented Moreau gradient with parameters $\stepsize, \augparam \geq 0$ as  
\begin{align*}
\nabla \anv_\augparam( f)(w; \stepsize) & = \argmin_{v \in \reals^\dimin} \stepsize f(w{-}v) {+} \frac{1}{2} \|v\|_2^2 {+} \frac{\augparam}{2} (f(w{-}v){-}f(w))^2
 = \argmin_{v \in \reals^\dimin} \augparam (f(w{-}v) {-} f(w) {+} \stepsize/\augparam)^2 {+}  \|v\|_2^2. 
\end{align*}	
For appropriate choices of $\stepsize, \augparam$, i.e., $\stepsize \leq \bar \stepsize/2$ and $\augparam \leq \bar \augparam /2$ with $\bar \augparam$ such that  $v\mapsto \bar \augparam (f(w-v)-f(w))^2/2 + \|v\|_2^2/2$ is convex, the minimizer is unique. 
As for the Moreau gradient, the augmented Moreau gradient defines a first-order optimality condition and an augmented Moreau gradient descent naturally ensures a decrease in the objective values. In practice, approximations of the augmented Moreau gradient may be obtained under appropriate conditions on $\stepsize, \augparam$, as detailed in Appendix~\ref{app:moreau}. The definition of augmented Moreau gradients sheds light on previous algorithms such as target propagation~\citep{lee2015difference} or proximal back-propagation~\citep{frerix2018proximal} while keeping the main properties of a classical Moreau envelope as detailed in Sec.~\ref{sec:rel_algos} and Appendix~\ref{app:moreau}.

\begin{figure}[t]
\begin{minipage}{0.5\linewidth}
	\begin{algorithm}[H]\caption{Forward pass \label{algo:gen_forward}}
		\begin{algorithmic}[1]
			\State {\bf Inputs:} Function $f$ parameterized by $(\dyn_t)_{t=1}^{\horizon}$ in~\eqref{eq:chain}, input $\state_0$, parameters $(\ctrl_t)_{t=1}^\horizon$
			\For{$t=1,\ldots, \horizon$}
			\State Compute $\state_t = \dyn_t(\ctrl_t, \state_{t-1})$
			\State Store $\state_{t-1}, \ctrl_t, \dyn_t$
			\EndFor 
			\State {\bf Output:}  Final result $\state_\horizon$
			\State {\bf Stored:} Intermediate comput. $(\state_{t-1}, \ctrl_t, \dyn_t)_{t=1}^\horizon$
		\end{algorithmic}
	\end{algorithm}
\end{minipage}~
\begin{minipage}{0.5\linewidth}
	\begin{algorithm}[H]\caption{Backward pass \label{algo:gen_backward}}
		\begin{algorithmic}[1]
			\State {\bf Inputs:} Stored $(\state_{t-1}, \ctrl_t, \dyn_t)_{t=1}^\horizon$, 
			last state $\state_\horizon$, objective $\obj$, stepsize $\stepsize$
			\State Initialize 
			$\lambda_\horizon =  \BP(\obj)(\state_\horizon, \stepsize)$ \label{line:obj_bp}
			\For{$t=\horizon, \ldots, 1$}
			\State 
			Compute $\lambda_{t-1} = \BP(\dyn_t( \ctrl_t, \cdot))(\state_t, \lambda_t)$ \label{line:state_bp}
			\State Compute $g_t = \BP(\dyn_t(\cdot, \state_t)) (\ctrl_t, \lambda_t)$ \label{line:ctrl_bp}
			\EndFor
			\State {\bf Output:} Oracle directions $(g_t)_{t=1}^\horizon$. 
		\end{algorithmic}
	\end{algorithm}
\end{minipage}
\end{figure}

\subsection{Differentiable Program for Moreau gradients}
We present first an overview of the proposed approximations of the Moreau gradients based on (i) generalizations of the notion of Moreau gradients to multivariate functions, (ii) back-propagation of these definitions along a graph of computation. 
Interpretations of the proposed procedure are presented in Sec.~\ref{sec:composition_moreau} and Sec.~\ref{sec:lag}.
\paragraph{Moreau gradient for multivariate functions}
For a multivariate function $f: \reals^ \dimin \rightarrow \reals^ \dimout$, a classical gradient encodes the linear form $\lambda \mapsto \nabla(\lambda^ \top f)(w)$. Similarly, we define the Moreau gradient and the augmented Moreau gradient of a multivariate function $f$ as the nonlinear forms
\[
\lambda \mapsto \nabla  \env(\lambda^ \top f)(w), \quad \lambda \mapsto \nabla \anv (f)(w; \lambda),
\]
where, in this case, $ \nabla \anv (f)(w; \lambda):= \argmin_{v \in \reals^\dimin}\left\{ \|f(w-v) - f(w) + \lambda\|^2 +  \|v\|_2^2 \right\}$.
Our goal in the following is to define a numerical program $\morprog$ which implements an approximation of the Moreau gradient or its augmented definition for a dynamical system such as~\eqref{eq:chain}. For example, for the Moreau gradient, the program $\morprog$ consists in 
\[
\morprog : \begin{cases}
\reals^\dimin & \rightarrow \reals^\dimout \times (\reals^\dimout \rightarrow \reals^ \dimin)\\
x & \mapsto (\diffprogfunc(w), \lambda \mapsto \widehat\nabla\env(\lambda^\top f)(w)),
\end{cases}
\]
which we aim to implement in a differentiable programming framework such that by evaluating $\texttt{out  = f(w)}$ we can access approximate Moreau gradients as $\widehat\nabla\env(\lambda^\top f)(w) = \texttt{automgrad(out, w, lambda)}$.

\paragraph{Back-propagation}
To implement Moreau gradients for functions of the form~\eqref{eq:chain}, we consider taking advantage of the structure of the problem just as a gradient oracle does by using automatic differentiation. A simplified overview of our approach is presented in Algo.~\ref{algo:gen_forward} and Algo.~\ref{algo:gen_backward}. Detailed implementations with additional hyper-parameters are presented in the following sections. 

The forward pass in Algo.~\ref{algo:gen_forward}, evaluates the function while keeping in memory the intermediate computations and the associated inputs as presented. During the backward pass in Algo.~\ref{algo:gen_backward}, we consider procedures that either use gradient-vector products, back-propagate the Moreau gradients or use regularized inversions of the intermediate computations. 
Namely, for a function $f \in \{\phi_t(w_t, \cdot), \phi_t(\cdot, x_t)\}$ evaluated at $z \in \{x_t, w_t\}$ respectively and a direction $\lambda$, we consider  back-propagation procedures $\BP$ of the form
\begin{align*}
	\GBP(f)(z, \lambda) & = \nabla f(z) \lambda\\
	\MBP(f)(z, \lambda) & = \widehat \nabla \env(\lambda^\top f)(z)  \approx \argmin_{y} \lambda^\top f(z-y) +\|y\|_2^2/2\\
	\IBP(f)(z, \lambda) & = \widehat \nabla \anv(f)(z; \lambda) \hspace{2pt} \approx \argmin_{y}  \|f(z- y) - f(z) + \lambda\|_2^2 +  \|y\|_2^2.
\end{align*}
Once oracle directions $(g_t)_{t=1}^\horizon$ are computed, the parameters $w_t$ can be updated using any optimization update as detailed in Sec.~\ref{sec:implem}. For an approximate Moreau gradient or augmented Moreau gradient descent, the variables are updated as
\begin{equation}\label{eq:updates}
	w_t  \leftarrow w_t -  g_t, \ \mbox{for} \ t \in \{1, \ldots, \horizon\}.
\end{equation}
The back-propagation rules $\MBP$, $\IBP$ are implemented with an optimization subroutine on their defining  problem as detailed in Sec.~\ref{sec:implem}. The rationale of the proposed approach can be understood by defining approximate chain rules for the Moreau gradient as explained in Sec.~\ref{sec:composition_moreau} or by considering Lagrangian or penalized formulations of the computation of the Moreau gradient as presented in Sec.~\ref{sec:lag}.

\begin{figure}[t]
	\begin{center}
\hspace*{50pt}\includegraphics[width=0.9\linewidth]{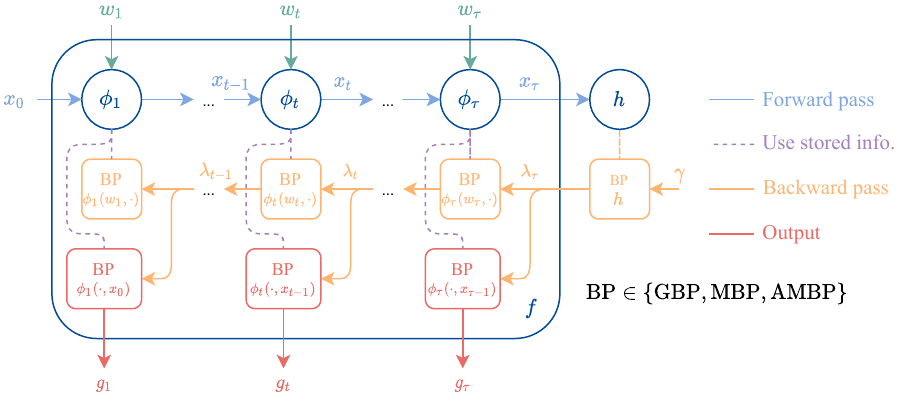}
		\caption{The proposed approximation of  Moreau gradients follows the computational scheme of gradient back-propagation, with   back-propagation rules $\MBP, \IBP$ replacing the usual gradient-vector product $\GBP$.\label{fig:soft_backprop}}
	\end{center}
\end{figure}

\begin{figure}[t]
	\begin{center}
		\includegraphics[width=0.45\linewidth]{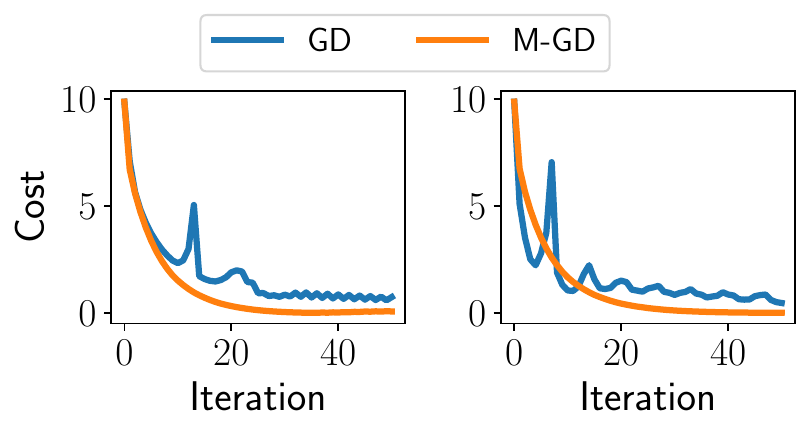}
		\vspace{-1ex}
		\caption{Gradient descent (GD) vs Moreau gradient descent (M-GD)  on the control of a pendulum.  \newline
			\textbf{Left:} horizon $\horizon=50$, \textbf{Right:} horizon $\horizon=100$\label{fig:ctrl}}
		\vspace{-2ex}
	\end{center}
\end{figure}

\begin{figure}[t]
	\begin{center}
	\raisebox{-0.5\height}{
		\includegraphics[width=0.5\linewidth]{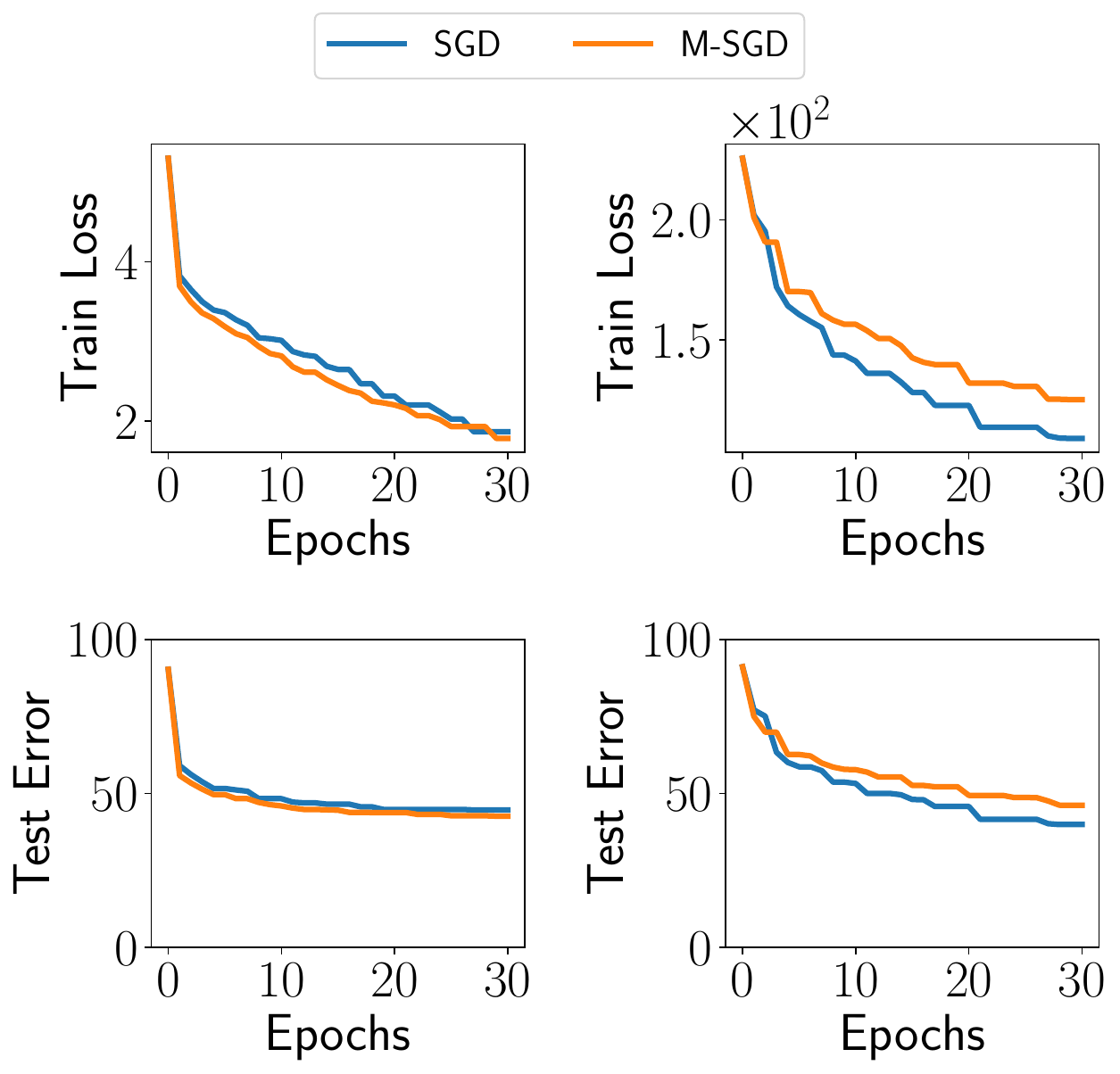}
	}~
	\raisebox{-0.5\height}{
		\includegraphics[width=0.33\linewidth]{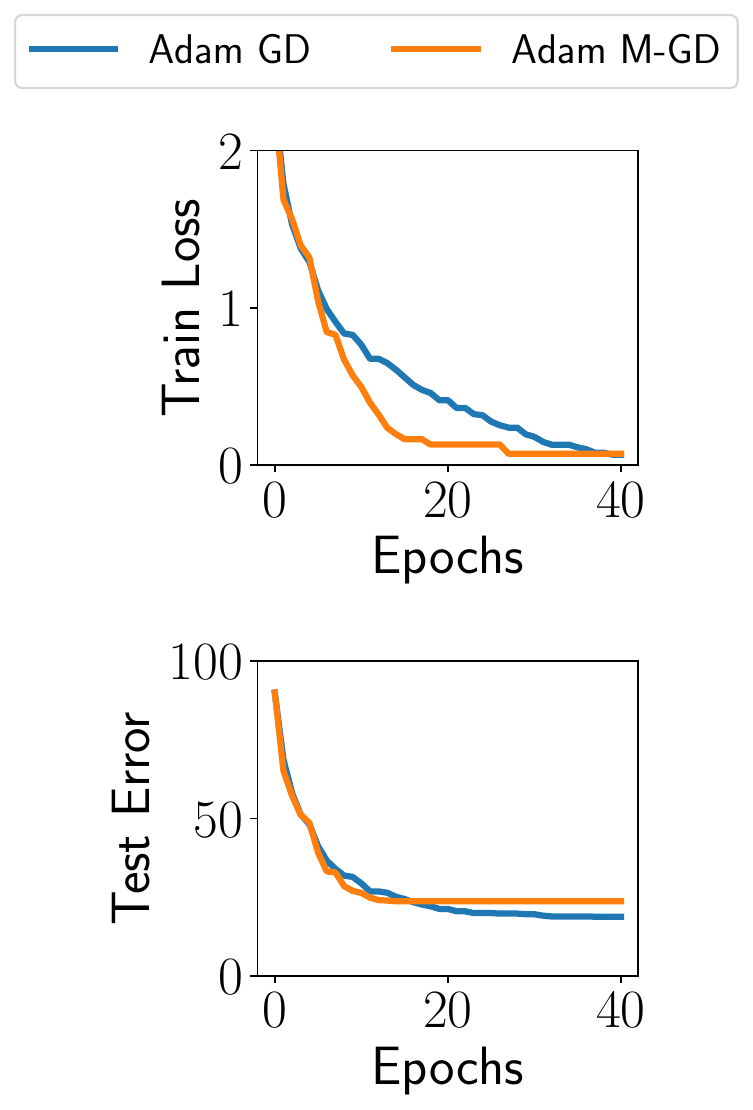}
	}
	\vspace{-1ex}
	\caption{\textbf{Left:} Stochastic Gradient Descent (SGD) versus Stochastic Moreau Gradient Descent  (M-SGD) on deep learning problems
		on CIFAR10 with an MLP,
		\textbf{Middle:} Stochastic Gradient Descent (SGD) versus Stochastic Moreau Gradient Descent  (M-SGD) on deep learning problems
		on CIFAR10 with a ConvNet as~\citet{frerix2018proximal},
		\textbf{Right:} Adam with Gradient oracle (Adam GD) versus Adam with Moreau gradients  (Adam M-GD) on CIFAR with AllCNN architecture. \label{fig:deep} \label{fig:proxprop}}
	\vspace{-2ex}
		\end{center}
\end{figure}

\paragraph{Numerical illustrations}
We illustrate the behavior of a Moreau gradient descent compared to gradient descent on (i) the control of a pendulum for various horizons $\tau$ in Fig.~\ref{fig:ctrl} and on (ii) an image classification task with various deep networks with stochastic surrogates of Moreau gradients and eventual additional momentum heuristics, as explained in Sec.~\ref{sec:oracle}, in Fig.~\ref{fig:deep}. We observe that the deterministic implementation of a Moreau gradient descent may provide smoother optimization for the control task, while the mini-batch stochastic counterpart of Moreau gradients compares favorably with stochastic gradient descent for the deep learning task. 
Experimental details on these illustrations are provided in Appendix~\ref{app:exp}. The rest of the paper focuses on presenting the rationale of the method and its implementation details. 

\section{Chain Rules}\label{sec:composition_moreau}
\subsection{Moreau Gradient}
\paragraph{Decomposition for a single composition}
We consider first the computation of the Moreau gradient of a single composition of $\infunc:\reals^\dimin \rightarrow \reals^\diminter$ and $\outfunc:\reals^\diminter\rightarrow \reals$. Under suitable assumptions on $\infunc$ and $\outfunc$, the following lemma presents a decomposition of the computation of the Moreau gradient similar to the computation of a classical gradient. See Appendix~\ref{ssec:app_simp_comp_moreau} for the proof and relaxed assumptions.

\begin{restatable}{lemma}{chainrulea}\label{prop:chain_rule}
	Consider $\infunc:\reals^\dimvar \rightarrow \reals^\diminter$  and $\outfunc:\reals^\diminter \rightarrow \reals$ Lipschitz-continuous and smooth with $\outfunc$ convex.  The Moreau gradient of $\outfunc\circ \infunc$ at a point $w$ exists for a stepsize $0 \leq \stepsize \leq 1/(2\lip_\outfunc \smooth_\infunc)$ and is given by
	\begin{align}
		\nabla \env(\stepsize \outfunc\circ \infunc)(w) &  =  \argmin_{v \in \reals^\dimvar}\left\{ {\adjvar^*}^\top \infunc(w  - v) + \frac{1}{2}\|v\|_2^2\right\} =  \nabla \env ({\adjvar^*}^ \top \infunc) (w),  \label{eq:dual_step}\\
		\mbox{where} \quad \adjvar^*  & = \argmax_{\adjvar \in \reals^\diminter} - (\stepsize \outfunc)^*(\adjvar) +  \env(\adjvar^\top \infunc) (w),
		\label{eq:dual_comput}
	\end{align}
	where $(\stepsize \outfunc)^ *$ is the convex conjugate of $\stepsize \outfunc$.
\end{restatable}
Compare~\eqref{eq:dual_step} to the classical gradient chain rule 
\begin{align*}
	\nabla (\stepsize \outfunc\circ \infunc)(w)  & =  \argmin_{v \in \reals^d}\left\{ \adjvar^\top \nabla  \infunc(w)^\top (w - v) + \frac{1}{2}\|v\|_2^2\right\}  \\
\mbox{where} \quad  \adjvar	& = \stepsize \nabla \outfunc( \infunc(w)).
\end{align*}
We retrieve the same structure, except that (i) for the Moreau gradient the dual direction $\adjvar^*$ is given by solving an optimization problem, (ii) the classical gradient minimizes a linear approximation of the inner function along this direction, while for the Moreau gradient the inner function itself is used. 
Note, in particular, that for $\outfunc, \infunc$ linear, the Moreau gradient matches the definition of the gradient of the composition. 

\paragraph{Chain rule approximation}
The computation of the dual direction $\adjvar^*$ in~\eqref{eq:dual_comput} amounts to solving 
\[
\max_{\adjvar \in \reals^\diminter} -(\stepsize h)^*(\adjvar) + e(\adjvar) \quad \mbox{for}\  e(\adjvar) = \env(\adjvar^ \top \infunc)(w).
\]
We consider approximating this problem by using  the proximal operator of $h$  which gives  access to the proximal operator of $\adjvar \rightarrow  -(\stepsize h)^*(\adjvar) $. Namely, we consider approximating the solution by a proximal gradient step from $0$, which, for some stepsize $\pgstepsize>0$, reads
\begin{align*}
	\adjvar^*  \approx \hat  \adjvar = \argmax_{\adjvar \in \reals^ k} \nabla e(0)^ \top \adjvar - ( \stepsize \outfunc)^ *(\adjvar) - \frac{1}{2\pgstepsize}\|\adjvar-0\|_2^ 2  = \pgstepsize \nabla \env(\pgstepsize^{-1}\stepsize \outfunc)(\infunc(w)),
\end{align*}
as detailed in Fact~\ref{fact:prox_grad}. Denoting $\scaling = \pgstepsize^{-1}\stepsize$, the resulting approximate Moreau gradient is then
\begin{align*}
	\nabla \env(\stepsize \outfunc\circ \infunc)(w)&  \approx \nabla \env (\stepsize {\hat \adjvar}^ \top \infunc) (w),  \quad
	\mbox{where} \quad \hat \adjvar =\scaling^{-1} \nabla \env(\scaling \outfunc)(\infunc(w)).
\end{align*}
More generally, for $\outfunc$ multivariate, we consider approximating the Moreau gradient as
\begin{align*}
	\nabla \env(\dir^\top \outfunc\circ \infunc)(w)&  \approx \nabla \env ({\hat \adjvar}^ \top \infunc) (w),  \quad
	\mbox{where} \quad \hat \adjvar  = \scaling^{-1}\nabla \env( \scaling \dir ^\top \outfunc)(\infunc(w)),
\end{align*}
that mimics the usual chain rule $\nabla (\dir^\top \outfunc\circ \infunc)(w) = \nabla (\adjvar^\top \infunc)(w)$ for $\adjvar = \nabla (\dir^\top \outfunc)(\infunc(w))$ up to a scaling parameter $\scaling$ by replacing $\nabla (\lambda^\top \infunc)$ by $\nabla \env(\lambda^\top f)$. The scaling parameter $\scaling$ is a hyperparameter that can be selected to ensure the feasibility of the computation of the intermediate Moreau gradients as detailed below for multiple compositions. While we consider only one step of a proximal gradient descent in order to build an approximation of the Moreau envelope for multiple compositions, we can also derive complexity bounds associated to the minimization of $\outfunc \circ \infunc$ using Moreau gradients as detailed in Appendix~\ref{app:chain_rule}.

\paragraph{Sequence of compositions}
The approach taken for a single composition amounts to computing the value of the inner function, i.e., $\infunc(w)$, then applying the Moreau gradient on the outer function $\outfunc$ on $f(w)$ to define a variable $\adjvar$ used to compute the Moreau gradient of the composition through the Moreau gradient of $f$. 
For multiple compositions of the form $\infunc = \varphi_\horizon \circ\ldots \circ  \varphi_1 : \reals^{d_0}  \rightarrow \reals^{d_\horizon}$ and $\obj: \reals^ {d_\horizon} \rightarrow \reals$, we consider generalizing this approach by (i) computing the intermediate values of the compositions $x_t = \varphi_t\circ \ldots  \varphi_1(x)$ for $t\in \{1, \ldots, \horizon\}$ in a forward pass, (ii) back-propagating some adjoint variables $\adjvar_t$ using the chain rule presented in the previous paragraph in a backward pass as
\begin{align}
	\adjvar_\horizon   & = \scaling_\horizon^{-1}\nabla \env(\scaling_\horizon h)(x_\horizon) \nonumber \\ 
 	\adjvar_{t-1}  &  = \scaling_{t-1}^ {-1} \nabla \env (\scaling_{t-1} \adjvar_t^ \top \varphi_t)(x_{t-1}) \label{eq:mor_grad_decomp} \ \mbox{for} \ t \in \{\horizon, \ldots, 2\} \nonumber  \\
	\nabla \env(\stepsize \obj \circ\infunc)(x)  & \approx \nabla \env(\stepsize \adjvar_1^ \top \varphi_1)(x). 
\end{align}
Each application of the proposed chain rule requires computing the Moreau gradient of some function with a given scaling parameter $\scaling_t$ chosen in advance as a hyperparameter. The scaling parameters can be chosen in such a way that the sub-problems defining the intermediate Moreau gradients are strongly convex, hence solvable with a first-order optimization method, as shown in Lemma~\ref{lem:stepsize}.
\begin{restatable}{lemma}{stepsizes}\label{lem:stepsize}
	Consider $\infunc = \varphi_\horizon \circ\ldots \circ  \varphi_1 : \reals^{d_0}  \rightarrow \reals^{d_\horizon}$, $\obj: \reals^ {d_\horizon} \rightarrow \reals$, and the back-propagation scheme presented in~\eqref{eq:mor_grad_decomp} for a given $x \in \reals^{d_0}$. Assume $h, \varphi_t$ to be Lipschitz continuous and smooth.
		For $t \in \{\horizon, \ldots, 1\}$, if $\scaling_t < 1/c_{t+1}$ with $c_t =  \smooth_{\varphi_t} \prod_{s=t+1}^ {\horizon+1} \lip_{\varphi_s}$ and $\varphi_{\horizon+1} = h$, the computation of $\adjvar_t$ in~\eqref{eq:mor_grad_decomp} amounts to a smooth strongly convex optimization problem solvable up to any accuracy by a first-order method at a logarithmic cost in the target accuracy. Similarly, if $\stepsize\leq c_1$, the computation of $\nabla \env(\stepsize \adjvar_1^ \top \varphi_1)(x)$ amounts also to a strongly convex optimization problem.
\end{restatable}
\begin{proof}
		For $t = \horizon$, the computation of $\adjvar_\horizon$ is a strongly convex optimization problem as long as $\scaling_\horizon < 1/\smooth_\outfunc$ as recalled in Lemma~\ref{lem:smooth_strg_cvx}.
		For $t\in\{1, \ldots, \horizon\}$, denote $x_t = \varphi_t \circ \ldots \circ \varphi_1(x) $ for a given $x \in \reals^{d_0}$.
		 For $t-1 \in \{\horizon-1, \ldots, 1\}$, the computation of the Moreau gradient $\nabla \env(\scaling_{t-1} \adjvar_t^ \top \varphi_t)$ needed for $\adjvar_{t-1}$ amounts to solving
	\[
	\min_{y_t} \ \scaling_{t-1} \adjvar_t^ \top \varphi_t(\state_{t-1} - y_{t-1}) + \frac{1}{2} \|y_{t-1}\|_2^ 2.
	\]
	The above problem is strongly convex if $\scaling_{t-1}< 1/\smooth_{\adjvar_t^\top \varphi_t} $, where $\smooth_{\adjvar_t^\top \varphi_t } $ is the smoothness constant of $\adjvar_t^\top \varphi_t$.
	We have 
	$
	\smooth_{\adjvar_t^\top \varphi_t } \leq \|\adjvar_t\|_2 \smooth_{\varphi_t},
	$
	so it remains to bound $\|\adjvar_t\|_2$. By definition of $\adjvar_t$, we have 
	\[
	\scaling_\horizon \nabla \outfunc(x_\horizon - \adjvar_\horizon) =  \scaling_\horizon \adjvar_\horizon, \quad  \scaling_{t-1}  \adjvar_t^\top \nabla \varphi_t(x_{t-1}- \adjvar_{t-1}) = \scaling_{t-1}\adjvar_{t-1}\quad \mbox{for} \ t \in \{\horizon, \ldots, 2\}.
	\]
	Hence for $t \in \{\horizon, \ldots, 1\}$, $\|\adjvar_t\|\leq \lip_\outfunc \prod_{s=t+1}^\horizon \lip_{\varphi_s}$. Hence $c_t$ defined in the statement of the lemma satisfies $c_t \geq \smooth_{\adjvar_t^\top \varphi_t } $ so taking $\scaling_{t-1} < 1/c_t $ ensures $\scaling_{t-1} <1 /\smooth_{\adjvar_t^\top \varphi_t}$ and the problem defining $\adjvar_{t-1}$ is strongly convex.  The claim for computing $ \nabla \env(\stepsize \adjvar_1^ \top \varphi_1)(x)$ follows from the same reasoning. 
\end{proof}

\paragraph{Dynamical systems}
Consider now  a dynamical system as in~\eqref{eq:chain}, i.e., a function $\infunc: \reals^\dimctrls \rightarrow \reals^{d_\horizon}$ parameterized by $\point \in \reals^{d_0}$ and $\horizon$  functions $\dyn_\layer:\reals^{\dimctrl_\layer} \times \reals^{d_{\layer-1}} \rightarrow \reals^{d_\layer}$ for $\layer=1,\ldots, \horizon$ such that for  $\ctrls = (\ctrl_1, \ldots, \ctrl_\horizon) \in \reals^\dimctrls, $ with $\ctrl_\layer \in \reals^{\dimctrl_\layer}$, the output of $\infunc$ is given by 
\begin{align}\label{eq:chain2}
	\infunc(\ctrls)  &= \state_\horizon,  \\
	\mbox{s.t.} \ \state_\layer & =	\dyn_\layer(\ctrl_\layer, \state_{\layer{-}1}) \ \mbox{for} \ \layer = 1, \ldots, \horizon.\nonumber
\end{align}
In this case, we consider computing the Moreau gradient of a composition $\outfunc \circ \infunc$ for $h:\reals^{d_\horizon} \rightarrow \reals$ with respect to every single parameter $w_t$. Namely, given a set of parameters $w\in \reals^p$, we consider computing the Moreau gradient of $\outfunc \circ \infunc_{w, t}$ for $t \in \{1, \ldots, \horizon\}$ where $f_{w,t}$ is defined by modifying the $t$\textsuperscript{th} parameter $w_t$ of the system while keeping all other parameters $(w_s)_{s\neq t}$ at their nominal values given by $w$, i.e,
\begin{align}
	\infunc_{w, t}(v_t) & = x_\horizon \label{eq:shift_coord}\\
		\mbox{s.t.} \  x_s & = \dyn_s(w_s, x_{s-1})  \quad \mbox{for} \ s \in \{t+1, \ldots, \horizon\} \nonumber\\
	x_t & = \dyn_t(v_t, x_{t-1}) \nonumber \\
	x_s & = \dyn_s(w_s, x_{s-1})  \quad \mbox{for} \ s \in \{1, \ldots, t-1\}. \nonumber
\end{align}
Denoting $(x_t)_{t=1}^\horizon$ the intermediate computations of~\eqref{eq:chain2} for a given set of parameters $w = (w_t)_{t=1}^\horizon$ and defining the functions $\varphi_s(y) = \phi_s(w_s, y)$ for $s \in \{1, \ldots \horizon\}$, the function $f_{w, t}$ of the dynamical system $f$ can be written as a composition
$
f_{w, t} = \varphi_\horizon \circ \ldots \varphi_{t+1} \circ \phi_t(\cdot, x_{t-1}).
$
The Moreau gradient of $\outfunc \circ \infunc_{w, t}$ on $w_t$ can then be computed with the  approach taken for multiple compositions as $\nabla \env (\stepsize_t \outfunc \circ \infunc_{w, t})(w_t) \approx g_t$ with
\begin{align}\label{eq:mbp_backprop_detailed}
	\adjvar_\horizon & =\scaling_\horizon^{-1} \nabla \env(\scaling_\horizon \outfunc) (x_\horizon) \nonumber\\
	\adjvar_{s-1} & = \scaling_{s-1}^{-1} \nabla \env(\scaling_{s-1} \adjvar_s^\top \varphi_s)(x_{s-1})  =  \scaling_{s-1}^{-1} \nabla \env(\scaling_{s-1} \lambda_s^\top \phi_s(w_s, \cdot))(x_{s-1})  \quad \mbox{for} \ s \in \{t+1, \ldots,  \horizon\} \nonumber\\
	g_t & = \nabla \env( \stepsize_t \lambda_t^\top \phi_t(\cdot, x_{t-1}))(w_t).
\end{align}
The computation of the intermediate adjoint variables $\adjvar_s$ are then shared when approximating $\nabla \env (\stepsize_t \outfunc \circ \infunc_{w, t})(w_t) $ for $t \in \{\horizon, \ldots, 1\}$. So to approximate $\nabla \env (\stepsize_t \outfunc \circ \infunc_{w, t})(w_t) $ for all $t \in \{\horizon, \ldots, 1\}$, the overall approach can be decomposed into (i) computing the intermediate values $(x_t)_{t=1}^\horizon$ of the dynamical system for $t\in \{1, \ldots, \horizon\}$ with the current parameters $w = (w_t)_{t=1}^\horizon$ in a \emph{forward pass}, (ii) back-propagate the adjoint variables $(\adjvar_t)_{t=1}^\horizon$ as in~\eqref{eq:mbp_backprop_detailed} in a \emph{backward pass}, (iii) output the approximate Moreau gradients $\nabla \env( \stepsize_t \lambda_t^\top \phi_t(\cdot, x_{t-1}))(w_t)$ for each $t\in \{1, \ldots, \horizon\}$ in the backward pass. We retrieve here the rationale presented in Algo.~\ref{algo:gen_forward} and Algo.~\ref{algo:gen_backward}. We detail the implementation of the backward pass used to approximate Moreau gradients in Algo.~\ref{algo:mbp_backward}. Details on the optimization subroutine, denoted $\mathcal{A}$ in Algo.~\ref{algo:mbp_backward}, are given in Sec.~\ref{sec:implem}. 

\begin{algorithm}\caption{Backward pass based on Moreau gradients \label{algo:mbp_backward}}
	\begin{algorithmic}[1]
		\State {\bf Inputs:} Stored $(\state_{t-1}, \ctrl_t, \dyn_t)_{t=1}^\horizon$, 
		last state $\state_\horizon$, objective $\obj$,  parameters $(\stepsize_t, \scaling_t)_{t=1}^\horizon$, optimization algorithm $\mathcal{A}$ such that $\mathcal{A}(f)$ is the result of an optimization algorithm applied to minimize $f$ and initialized at $0$
		\State {\bf Notation:} For a function $f$, $\widehat \nabla \env(f)(x) = \mathcal{A}(f(x-\cdot) + \|\cdot\|_2^2/2$) 
		\State Initialize $\adjvar_\horizon = \scaling_\horizon^{-1} \widehat \nabla \env(\scaling_\horizon \obj)(\state_\horizon)$
		\For{$t=\horizon, \ldots 1$}
		\State Compute $\adjvar_{t-1} = \scaling_{t-1}^{-1} \widehat \nabla \env(\scaling_{t-1} \adjvar_t^\top \phi_t(w_t, \cdot))(x_{t-1})$
		\State Compute $g_t = \widehat \nabla \env(\stepsize_t \adjvar_t^\top \phi_t(\cdot, x_{t-1}))(w_t)$ \label{line:oracle_mbp}
		\EndFor 
		\State {\bf Output:} Oracle directions $g_1, \ldots, g_\horizon$ 
	\end{algorithmic}
\end{algorithm}

\subsection{Augmented Moreau Gradients}
In this section, we present a simple chain rule to approximate the simplified augmented Moreau gradient presented in Sec.~\ref{sec:moreau}. A refined implementation based on a penalized formulation of the problem is presented in Sec.~\ref{sec:lag} with detailed pseudo-code. 

\paragraph{Decomposition for a single composition}
The computation of the augmented Moreau gradient $\nabla \anv(\outfunc\circ \infunc)(w; \stepsize) $ of the composition $\outfunc\circ \infunc$ can be decomposed as
\begin{align} \label{eq:amg_decomp}
\min_{u \in \reals^\diminter} \ & \left(\outfunc(\infunc(w) - u_1)  - \outfunc(\infunc(w)) + \stepsize\right)^2/2 +  p(u_2), \quad \mbox{where} \ p(u) = \min_{v \in \reals^\diminter  } \|v\|_2^2/2.\\
\mbox{s.t.} \hspace{5pt} & u_1 = u_2  \hspace{228pt} \mbox{s.t.} \  \infunc(w)  - \infunc(w -v) = u. \nonumber
\end{align}
The above formulation is then akin to the computation of $\nabla \anv(\outfunc)(\infunc(w); \stepsize)$ up to a modified regularization and can be tackled by methods for sums of functions such as the Alternate Direction Method of Multipliers (ADMM)~\citep{boyd2011distributed} as explained below. 

\paragraph{Chain rule approximation}
For an objective of the form $\min_x f_1(x) + f_2(x)$, the ADMM method proceeds by considering the constrained problem $\min_{x_1 =x_2} f_1(x_1) + f_2(x_2)$ and uses a primal-dual method on the augmented Lagrangian, $\lag_{\aug, \rho}(x_1, x_2, \lambda) = f_1(x_1) + f_2(x_2) + \lambda^\top(x_1-x_2) + \rho\|x_1-x_2\|_2^2/2$ with the primal method defined by alternate minimization on the augmented Lagrangian,  i.e., $x_1 \leftarrow \argmin_{x} f_1(x) + \lambda^\top(x-x_2) + \rho\|x-x_2\|_2^2/2 $, $x_2 \leftarrow \argmin_{x} f_2(x) + \lambda^\top(x_1-x) + \rho\|x_1-x\|_2^2/2$ and $\lambda \leftarrow \lambda + \rho(x_1-x_2)$, where $\lambda$ encapsulates Lagrange multipliers, i.e., dual variables~\citep{boyd2011distributed}. For problem~\eqref{eq:amg_decomp}, one iteration of the ADMM approach starting from $u_1 = u_2 = \lambda = 0$ with $\rho=1$ gives
\begin{align*}
	u_1 &= \argmin_{u \in \reals^\diminter}\left\{\left(\outfunc(\infunc(w) - u)  - \outfunc(\infunc(w)) + \stepsize\right)^2 + \|u\|_2^2\right\} \\
	u_2 & = \argmin_{u \in \reals^\diminter}  \left\{p(u) + \|u-u_1\|_2^2/2 \right\} = \infunc(w)  - \infunc(w- \hat v), \ \mbox{for}\ 
	\hat v = \argmin_{v\in \reals^\dimin} \left\{\|\infunc(w) -\infunc(w {-}v) - u_1\|_2^2 + \|v\|_2^2   \right\}.
\end{align*}
We ignore the update on the dual variable here and simply consider  $\hat v$ as an approximation of the augmented Moreau gradient whose computation can be summarized as
\begin{align*}
		\nabla \env(\outfunc\circ \infunc)(w;\stepsize)   \approx \hat v = \nabla \anv(\infunc) (w; \lambda), \quad
\mbox{where} \quad	\lambda =\nabla \anv(\outfunc) (\infunc(w); \stepsize).
\end{align*}
More generally for $h$ multivariate, we consider approximating the augmented Moreau gradient as
\begin{align*}
	\nabla \env(\mu^\top \outfunc\circ \infunc)(w; \lambda) & \approx  \nabla \anv(\infunc) (w;  \lambda), \hspace{16pt}
	\mbox{where} \hspace{10pt} \lambda  = \nabla \anv(\outfunc) (\infunc(w);  \mu),
\end{align*}
which  mimics the usual chain rule $\nabla (\mu^\top \outfunc \circ \infunc)(w) = \nabla (\lambda^\top f)(w)  $ for $\lambda = \nabla (\mu^\top h)(f(w))$ by replacing the operators $\nabla (\lambda^\top \infunc)$ with $\nabla \anv(\infunc)(\cdot;\lambda)$.

\paragraph{Multiple compositions and dynamical systems}
For a sequence of compositions $\infunc = \varphi_\horizon \circ\ldots \circ  \varphi_1 : \reals^{d_0}  \rightarrow \reals^{d_\horizon}$ and $\obj: \reals^ {d_\horizon} \rightarrow \reals$, we can apply the proposed chain rule iteratively to approximate $\nabla \anv(h\circ f)(x;\stepsize)$ as $\nabla \anv( \obj \circ\infunc)(x;\stepsize)  \approx \lambda_0$ with
\begin{align}\label{eq:aug_morgrad_seq_comp}
	\adjvar_\horizon   & = \nabla \anv(h)(x_\horizon; \stepsize), \; 
	\adjvar_t   = \nabla \anv (\varphi_{t+1})(x_t; \adjvar_{t+1}) \ \mbox{for} \ t \in \{\horizon-1, \ldots, 0\}, 
\end{align}
for $x_1, \ldots, x_\horizon$ defined as $x_t = \varphi_t(x_{t-1})$  for $t=0, \ldots, \horizon-1$ with $x_0 = x$.  

For dynamical systems as in~\eqref{eq:chain2}, we consider computing $\nabla \anv ( \outfunc \circ \infunc_{w, t})(w_t;\stepsize_t)$ for $t\in \{1, \ldots, \horizon\}$ and $f_{w, t}$ defined as in~\eqref{eq:shift_coord} by applying the rule~\eqref{eq:aug_morgrad_seq_comp} on $
f_{w, t} = \varphi_\horizon \circ \ldots \varphi_{t+1} \circ \phi_t(\cdot, x_{t-1})
$ and $\varphi_s(y) = \phi_s(w_s, y)$ for $s \in \{1, \ldots \horizon\}$. Such an approach is summarized by Algo~\ref{algo:gen_forward} and Algo.~\ref{algo:gen_backward} with the procedure $\IBP$.

\section{Lagrangian and Penalized Viewpoints}\label{sec:lag}
The classical gradient back-propagation can be interpreted from numerous viewpoints such as computing the Lagrange multipliers associated with a linearization of the objective~\citep[Section 2.6]{bertsekas1999nonlinear} or making a coordinate gradient descent on a penalized formulation of the dynamical constraints~\citep[Proposition 1]{frerix2018proximal}. Following these approaches, we interpret the proposed approximation of the Moreau gradients through the lens of a Lagrangian or a penalized formulation of the dynamical constraints. 

\subsection{Lagrangian Viewpoint}
Consider the minimization of $\outfunc \circ \infunc$ with $\infunc$ defined as in~\eqref{eq:chain2} by $\horizon$ functions $(\phi_t)_{t=1}^\horizon$ such that, for $w = (w_1, \ldots, w_\horizon)$, $f(w) = x_\horizon$ with $x_t = \phi_t(w_t, x_{t-1})$ for $t \in \{1, \ldots, \horizon\}$ and $x_0$ fixed. The minimization of $\outfunc \circ \infunc$ amounts to a saddle point problem defined by a Lagrangian $\mathcal{L}$ as follows
\begin{equation}\label{eq:lag}
	\min_{w = (w_1, \ldots, w_\horizon)} \ \outfunc(\infunc(w)) = 
	\min_{\substack{w_1, \ldots, w_\horizon \\ 
			x_1, \ldots, x_\horizon}}  \sup_{\adjvar_1, \ldots, \adjvar_\horizon }\left\{\mathcal{L}(w_{1:\horizon}, x_{1:\horizon}, \adjvar_{1:\horizon}) = \outfunc(x_\horizon) + \sum_{t=1}^{\horizon} \adjvar_t^\top (\phi_t(w_t, x_{t-1}) - x_t)  \right\},
\end{equation}
for $x_0$ fixed.  An approximate Moreau gradient step on $\outfunc \circ \infunc$ from $w=(w_1, \ldots, w_\horizon)$, i.e., $w_t \leftarrow w_t - g_t$ with  $g_t$ defined as in~\eqref{eq:mbp_backprop_detailed} amounts to a block coordinate Moreau gradient descent on the Lagrangian as formally stated and explained in Lemma~\ref{lem:lag}. In the vocabulary of proximal operators, a block coordinate Moreau gradient descent on the Lagrangian amounts to a block coordinate proximal point method on the Lagrangian. 

\begin{lemma}\label{lem:lag}
	For $h:\reals^{d_\horizon} \rightarrow \reals$ and $f:\reals^p \rightarrow \reals^{d_\horizon}$
	defined as in~\eqref{eq:chain2},
	consider an  approximate Moreau gradient step on $\outfunc \circ \infunc$ from $w=(w_1, \ldots, w_\horizon)$, defined as $w_t \leftarrow w_t - g_t$ for $g_t$ defined as in~\eqref{eq:mbp_backprop_detailed} and $t \in \{1, \ldots, \horizon\}$.  This step amounts to 
	a block coordinate Moreau gradient pass on the Lagrangian~\eqref{eq:lag} starting from $w_{1:\horizon} = (w_1, \ldots, w_\horizon)$, $\lambda_{1:\horizon} = (0, \ldots, 0)$ and $x_{1:\horizon}$ defined by $x_t = \phi_t(w_t, x_{t-1})$ for $t \in \{1, \ldots, \horizon\}$.
\end{lemma}
\begin{proof}
	Denote by a superscript ${}^+$ the updated parameters along the block coordinate Moreau gradient pass applied to the parameters of the Lagrangian~\eqref{eq:lag}. We consider a block coordinate Moreau gradient pass whose updates are $x_\horizon \rightarrow \lambda_\horizon \rightarrow (w_\horizon, x_{\horizon-1}) \rightarrow \lambda_{\horizon-1} \rightarrow \ldots \rightarrow (w_1, x_0) $ where $\rightarrow$ denoted the next update to be done and $(w_t, x_{t-1})$ mean that the updates are performed simultaneously by a block-coordinate inner scheme.
	At layer $\horizon$, a Moreau gradient step on $\state_\horizon$ with step-size $\auxstepsize_\horizon$,  is given as, using that $\costate_\horizon=0$ and $\state_\horizon= \dyn_\horizon(\ctrl_\horizon, \state_{\horizon-1})$,
	\begin{align}
		\state_\horizon^+ 
		& = \state_\horizon - \nabla \env(\auxstepsize_\horizon \obj)(\state_\horizon). 
		\label{eq:lag_state_final}
	\end{align}
	At each layer $t\in \{\horizon,\ldots, 1\}$,
	\begin{enumerate}[nosep, leftmargin=12pt]
		\item a Moreau gradient step on $\costate_t=0$ with step-size $\auxxstepsize_t$ is given as, for $\state_t^ +, \state_{t-1}, \ctrl_t$ fixed,  denoting $\ell_t(\lambda) = \lambda^\top(\dyn_t(\ctrl_t,  \state_{t-1}) - \state_{t}^ +)$, 
		\begin{align}
			\costate_t^+ &  = \costate_t + \auxxstepsize_t \nabla \env(\ell_t)(\lambda_t) = 0 + \auxxstepsize_t (\dyn_t(\ctrl_t,  \state_{t-1}) - \state_{t}^ +) = \auxxstepsize_t(\state_t -  \state_t^ +),
			\label{eq:lag_costate}
		\end{align}
		using that $\state_t = \dyn_t( \ctrl_t,  \state_{t-1})$ and that $\nabla \env(\ell) = \nabla \ell$ for $\ell$ a linear function,
		\item a Moreau gradient step on $\ctrl_t$ with step-size $\stepsize_t$  reads, for $\costate_t^ +, \state_{t-1}$ fixed, 
		\begin{align*} 
			\ctrl_t^+ 
			& = \ctrl_t - \nabla \env({\stepsize_t \lambda_t^ +}^ \top \phi_t(\cdot, \state_{t-1}) )(\ctrl_t),
		\end{align*}
		
		\item  if $t\in \{\horizon \ldots, 2\}$, a  Moreau gradient step  on $\state_{t-1}$ with step-size $\auxstepsize_{t-1}$ reads, for $\costate_t^ +, \ctrl_t, \costate_{t-1} = 0$ fixed,
		\begin{align}
			\state_{t-1}^+ 
			& = \state_{t-1} -  \nabla \env(\auxstepsize_{t-1}{\lambda_t^ +}^ \top \phi_t(\ctrl_t, \cdot))(\state_{t-1}).
			\label{eq:lag_state}
		\end{align}
		
	\end{enumerate}
	By combining~\eqref{eq:lag_state_final},~\eqref{eq:lag_costate} and~\eqref{eq:lag_state}, the above decomposition can be rewritten as $\lambda_\horizon^+ = \auxxstepsize_\horizon \nabla \env(\auxstepsize_\horizon \obj)(\state_\horizon)$ followed by $\lambda_{t-1}^+ = \auxxstepsize_t \nabla \env (\auxstepsize_{t-1}{\lambda_t^ +}^ \top \phi_t(\ctrl_t, \cdot))(\state_{t-1})$ for $t=\horizon, \ldots, 1$. By taking $\auxxstepsize_t = \auxstepsize_t^{-1}$, we retrieve the decomposition~\eqref{eq:mbp_backprop_detailed}. 
\end{proof}
The Lagrangian viewpoint justifies the form of the updates taken by approximate Moreau gradient steps as shown in Lemma~\ref{lem:optmor}. 
\begin{restatable}{lemma}{optmor}\label{lem:optmor}
	For $h:\reals^{d_\horizon} \rightarrow \reals$ and $f:\reals^p \rightarrow \reals^{d_\horizon}$
	defined as in~\eqref{eq:chain2},
	 consider the computation of $\nabla\env (\stepsize_t\outfunc \circ \infunc_{w, t})(w_t)$ with $\infunc_{w, t}$ defined in~\eqref{eq:shift_coord} for a given $w \in \reals^p$. If $\stepsize_t \leq 1/c_t$ with $c_t = \lip_h \smooth_{\phi_t} \prod_{s=t+1}^ \horizon \lip_{\phi_s}$, the computation of $\nabla\env (\stepsize_t\outfunc \circ \infunc_{w, t})(w_t)$ is given by 
	\begin{equation}\label{eq:mor_grad_opt}
		\nabla\env (\stepsize_t\outfunc \circ \infunc_{w, t})(w_t)= \argmin_{\auxctrl_t \in \reals^{\dimctrl_t}} \  {\costate_t^*}^\top \dyn_t( \ctrl_t - \auxctrl_t, \state_{t-1})  {+} \|\auxctrl_t\|_2^2/(2\stepsize_t),
	\end{equation}	
	for $\state_{t-1}$ defined by $x_s = \phi_t(\ctrl_t, x_{s-1})$ for $s \in \{1, \ldots, t-1\}$, with $\lambda_{t:\horizon}^ *$ and $x_{t:\horizon}^ *$ regular solutions of 
	\begin{align}\label{eq:lag_opt}
		\min_{\auxctrl_t }\min_{x_{t:\horizon}} \sup_{\costate_{t:\horizon}} &\   h(x_\horizon) + \sum_{s=t+1}^{\horizon} \lambda_s^\top(\phi_s(\ctrl_s, x_{s-1}) - x_s)   + \lambda_t^ \top(\phi_t(\ctrl_t + \auxctrl_t, \state_{t-1}) - x_t) + \frac{1}{2\stepsize}\|\auxctrl_t\|_2^ 2.
	\end{align}
\end{restatable}
\begin{proof}
	If $\auxctrl_t^*, \costates_{t:\horizon}^*, \states_{t:\horizon}^*$ are regular solutions of~\eqref{eq:lag_opt}, they satisfy the associated Karush-Kuhn-Tucker conditions. In particular, we have for $s\in \{t+1,\ldots, \horizon\}$,
	\begin{align}\label{eq:kkt}
		\costate_\horizon^*  = \nabla \obj(\state_\horizon^*), \quad 
		\costate_{s-1}^* =  \nabla_{\state_{s-1}^*}  \dyn_s(\ctrl_s, \state_{s-1}^*) \costate_s^*, \quad
		\stepsize_t\nabla_{\ctrl_t-v_t^*}  \dyn_t(\ctrl_t - \auxctrl_t^*, \state_{t-1}) \costate_t^*  =  \auxctrl_t^*. 
	\end{align}
	We then have  for $s\in \{t,\ldots, \horizon\}$, $\|\costate_s^*\|_2 \leq \lip_\obj\prod_{j=s+1}^\horizon \lip_{\dyn_j}$.
	Therefore if $\stepsize_t<c_t^{-1}$, the function $\auxctrl_t \rightarrow {\lambda_t^*}^\top\dyn_t(\ctrl_t - \auxctrl_t) + \|\auxctrl_t\|_2^2/(2\stepsize_t)$ is strongly convex and $\auxctrl_t^*$ is also given as the solution of~\eqref{eq:mor_grad_opt}.	
\end{proof}

The above proposition emphasizes that the dual directions $\lambda_t^ *$ that define the Moreau gradient $\nabla \env(\stepsize_t \obj \circ f_{w, t})(w_t)  {=} \nabla \env(\stepsize_t {\lambda_t^ *}^ \top \phi_t(\cdot, \state_{t-1}))(\ctrl_t)$ are a priori given as the solution of an optimization problem for each $t$. In our implementation, we consider approximating $\lambda_t$ by means of a backward  pass as defined in Algo.~\ref{algo:mbp_backward} that takes advantage of the structure of the problem. 

\subsection{Penalized Viewpoint}
Consider again the minimization of $\outfunc \circ \infunc$ with $\infunc$ defined as in~\eqref{eq:chain2}.
Rather than considering the Lagrangian as in~\eqref{eq:lag}, we can consider a penalized formulation of the problem as defined below
\begin{equation}\label{eq:pen}
	\pen_{\reglag}(w_{1:\horizon}, x_{1:\horizon}) =  h(x_\horizon) + \sum_{t=1}^{\horizon} \frac{\reglag}{2}\|\dyn_t(w_t, x_{t-1}) -  x_t\|_2^2,
\end{equation}
with $x_0$ fixed and $\reglag> 0$ a penalty parameter. We can then adapt the approach presented in Lemma~\ref{lem:lag} to the penalized formulation by considering
a block coordinate Moreau gradient pass on the penalized formulation. Such an approach can be expressed in terms of the Moreau gradient of the objective $h$ and the augmented Moreau gradients of the intermediate functions as shown in the following lemma. 

\begin{lemma}\label{lem:pen}
For $h:\reals^{d_\horizon} \rightarrow \reals$ and $f:\reals^p \rightarrow \reals^{d_\horizon}$ defined as in~\eqref{eq:chain2},
a block coordinate Moreau gradient pass on the penalized objective~\eqref{eq:pen} starting from $w_{1:\horizon} = (w_1, \ldots, w_\horizon)$,  and $x_{1:\horizon}$ defined by $x_t = \phi_t(w_t, x_{t-1})$ for $t \in \{1, \ldots, \horizon\}$ amounts to computing
	\begin{align*}
	\adjvar_\horizon & = \nabla\env(\scaling_\horizon \outfunc)(x_\horizon), \quad \adjvar_{t-1} = \nabla \anv_{\scaling_{t-1}\reglag} (\phi_t(w_t, \cdot))(x_{t-1}; \scaling_{t-1}\reglag \lambda_t) \ \mbox{for} \ t \in \{2, \ldots, \horizon\} \\
	w_t ^+& = w_t - g_t \hspace{62pt} g_t =  \nabla \anv_{\stepsize_t\reglag}(\phi_t(\cdot, x_{t-1}))(w_t; \stepsize_t\reglag\lambda_t) \hspace{10pt} \mbox{for} \ t \in \{1,\ldots, \horizon\},
	\end{align*}
where $\scaling_t, \stepsize_t \geq 0$ for $t\in \{1, \ldots, \horizon\}$ are some hyper-parameters.
\end{lemma}
\begin{proof}
	Denote by a superscript ${}^+$ the updated parameters along the block coordinate Moreau gradient pass applied to the parameters of the penalized  formulation defined in~\eqref{eq:pen}.  We consider a block coordinate Moreau gradient pass whose updates are $x_\horizon \rightarrow  (w_\horizon, x_{\horizon-1})   \rightarrow \ldots \rightarrow (w_1, x_0) $ where $\rightarrow$ denoted the next update to be done and $(w_t, x_{t-1})$ mean that the updates are performed simultaneously in an inner block coordinate scheme. 
	At layer $\horizon$, a Moreau gradient step on $\state_\horizon$ with step-size $\auxstepsize_\horizon$,  is given as, using that $\dyn_\horizon(\ctrl_\horizon, \state_{\horizon-1}) = \state_\horizon$ for $w_\horizon, x_{\horizon-1}$ fixed,
	\begin{align*}
		\state_\horizon^+ & = \state_\horizon - \nabla \env(\auxstepsize_\horizon (h+ \reglag \|\cdot-x_\horizon\|^2/2))(x_\horizon) = \state_\horizon - \nabla \env(\scaling_\horizon h)(x_\horizon),
	\end{align*}
where $\scaling_\horizon = (\reglag+ \auxstepsize_\horizon^{-1})^{-1}$.
	At each layer $t\in \{\horizon,\ldots, 1\}$,
	\begin{enumerate}[nosep, leftmargin=12pt]
		\item a Moreau gradient  step on $\ctrl_t$ with step-size $\stepsize_t$  reads, for $\state_t^+, \state_{t-1}$ fixed,
		\begin{align*} 
			\ctrl_t^+ & = \ctrl_t - \ \argmin_{\auxctrl_t } \ \frac{\stepsize_t\reglag}{2}\|\dyn_t(\ctrl_t -\auxctrl_t,  \state_{t-1}) -  \state_t^ {+}\|_2^2 +  \frac{1}{2}\|\auxctrl_t \|_2^2 = w_t - \nabla \anv_{\stepsize_t\reglag}(\phi_t(\cdot, x_{t-1}))(w_t;\reglag \stepsize_t\lambda_t),
		\end{align*}
		where we denoted $\lambda_t = x_t- x_t^+$,
		\item  if $t\in \{\horizon, \ldots, 2\}$, a Moreau gradient step on $\state_{t-1}$ with step-size $\auxstepsize_{t-1}$ for $\state_t^+, \ctrl_t , \state_{t-2}, \ctrl_{t-1}$ fixed such that  $\state_{t-1}= \dyn_{t-1}(\ctrl_{t-1}, \state_{t-2})$, reads 
		\begin{align*}
			\state_{t-1}^+  & = \state_{t-1} -  \argmin_{\auxstate_{t-1} } \frac{1 + \auxstepsize_{t-1}\reglag}{2 }\|\auxstate_{t-1}\|^2 + \frac{\auxstepsize_{t-1}\reglag}{2}\|\dyn_t( \ctrl_t,\state_{t-1} -\auxstate_{t-1}) -  \state_t^ {+} \|_2^2 \\
			& =  \state_{t-1} -  \nabla \anv_{\scaling_{t-1} \reglag} (\phi_t(w_t, \cdot))(x_{t-1}; \reglag\scaling_{t-1} \lambda_t),
		\end{align*}
	where $\scaling_{t-1}  = (\reglag+ \auxstepsize_{t-1}^{-1})^{-1}$ and $\lambda_t = x_t- x_t^+$. 
		By expressing the pass in terms of the variables $\lambda_t$ defined above we get the claimed procedure.
	\end{enumerate}
\end{proof}
An oracle based on such a block coordinate Moreau gradient pass on the penalized formulation records the inputs and the programs used in a forward pass as in Algo.~\ref{algo:gen_forward} and back-propagate adjoint variables using the augmented Moreau gradients of the intermediate computations as presented in detail in~Algo.~\ref{algo:ambp_backward}. We retrieve the same structure as in Algo.~\ref{algo:mbp_backward}, i.e., recursive applications of non-linear operators based on the definition of a Moreau envelope with some additional hyper-parameters that can be used to ensure that the intermediate problems are strongly convex. One difference is that here the first oracle used is a Moreau gradient on the outer function $\outfunc$ rather than an augmented Moreau gradient.

\begin{algorithm}\caption{Backward pass based on augmented Moreau gradients \label{algo:ambp_backward}}
	\begin{algorithmic}
		\State {\bf Inputs:} Stored $(\state_{t-1}, \ctrl_t, \dyn_t)_{t=1}^\horizon$, 
last state $\state_\horizon$, objective $\obj$,  parameters $(\stepsize_t, \scaling_t)_{t=1}^\horizon$, optimization algorithm $\mathcal{A}$ such that $\mathcal{A}(f)$ is the result of an optimization algorithm applied to minimize $f$ and initialized at $0$
\State {\bf Notation:} For a function $f$, 

$\widehat \nabla \env(f)(x) = \mathcal{A}(f(x-\cdot) + \|\cdot\|_2^2/2)$, $\widehat\nabla \anv_\reglag (f)(x; \lambda)  = \mathcal{A}(\reglag\|f(x-\cdot) - f(x) + \lambda/\reglag\|_2^2 + \|\cdot\|_2^2)$
\State Initialize $\adjvar_\horizon = \widehat \nabla\env(\scaling_\horizon \outfunc)(x_\horizon)$
\For{$t=\horizon, \ldots 1$}
\State Compute 
$\adjvar_{t-1} =\widehat \nabla \anv_{\scaling_{t-1}\reglag} (\phi_t(w_t, \cdot))(x_{t-1}; \scaling_{t-1}\reglag\lambda_t) $
\State Compute $g_t = \widehat\nabla \anv_{\stepsize_t\reglag}(\phi_t(\cdot, x_{t-1}))(w_t;\stepsize_t\reglag \lambda_t )$
\EndFor 
\State {\bf Output:} Oracle directions $g_1, \ldots, g_\horizon$
	\end{algorithmic}
\end{algorithm}

\section{Implementation Details}\label{sec:implem}
The implementation of the proposed oracles as presented in detail in, e.g., Algo.~\ref{algo:mbp_backward} requires having access to an approximation of the Moreau gradients of the intermediate computations in closed form  or by means of an optimization subroutine. This additional flexibility can be exploited to develop various computational schemes and to control the computational complexity of the proposed approach as presented below.

\subsection{Moreau Gradients in Closed Form}
If the Moreau gradients of the intermediate computations can be computed in closed form, the overall computational complexities of our approximations of the Moreau gradient are the same as the ones of a gradient back-propagation algorithms, up to the overhead of the closed-form computation. We present several cases where the Moreau gradients can be computed in closed-form below. 

\paragraph{Linear case} The Moreau gradients of linear functions or more generally affine functions such as $\ell: w \rightarrow Aw +b$ coincide with the usual gradients in the sense that $\nabla \env(\lambda^\top \ell)(w) = \argmin_v \lambda^\top A(w-v) + \|v\|_2^2/2 = A^\top \lambda = \nabla (\lambda^\top \ell)(w)$. 
 Similarly,  the augmented Moreau gradients can be computed in closed form as we have $\nabla \anv_\reglag(\ell)(w; \lambda) = \argmin_v \{\reglag\|Av-\lambda/\reglag\|_2^2 + \|v\|_2^2\} = (\reglag A^\top A + I)^{-1} A^\top \lambda$. For the augmented Moreau gradients the computation requires the inversion of a matrix at a computational cost of $O(n^3)$, while the Moreau gradient, i.e., the usual matrix-vector product costs at most $O(n^2)$. 

\paragraph{Nonlinear simple functions}
More generally, for nonlinear functions that are the element-wise application of a simple operation, we can have access to the Moreau gradients in closed form. Consider for example the ReLU function $\relu: w \mapsto \max\{w, 0\}$ and its element-wise application $r: w \mapsto (\relu(w_i))_{i=1}^n$. The Moreau gradient can be computed for any $\stepsize$ as $\nabla\env(\stepsize \relu)(x) =  \min\{\relu(x), \relu(\stepsize)\}) + \min\{\stepsize, 0\} \ones_{x\geq \stepsize/2}(x)$. Similarly the augmented Moreau gradient is given as $\nabla \anv_\reglag(\relu)(x; \stepsize)  = \min\{\relu(x), \relu(\stepsize/(\reglag +1))\} + \min\{\stepsize/(\reglag +1) , 0\} \ones_{x\geq 0}(x)$. The element-wise application of the ReLU has then a Moreau gradient given  for any $\lambda \in \reals^d$ as $\nabla \env(\lambda^\top r)(w) = ( \min\{\relu(w_i), \relu(\lambda_i)\}) + \min\{\lambda_i, 0\} \ones_{w_i\geq \lambda_i/2} )_{i=1}^d$. Similarly, the Moreau gradients of tangent hyperbolic activations or sigmoid activations can be found by analyzing the corresponding univariate functions.

\subsection{Moreau Gradients Approximation Schemes}
While we may have access to the Moreau gradients of some elementary functions in closed form, our framework also adapts to any function by considering a subroutine able to approximate the Moreau gradients. This offers flexibility on the scale at which we consider approximating the Moreau gradients. Consider for example $\phi(x) = r(Ax+b)$ with $r$ the element-wise application of the ReLU activation function. Both the ReLU activation and the affine function $x\rightarrow Ax+b$ admit closed-form expressions for their Moreau gradients so that we can apply the proposed chain rule presented in Sec.~\ref{sec:composition_moreau}. On the other hand, we can consider approximating the Moreau gradient of $\phi$ directly by means of an optimization subroutine as $\nabla \env(\lambda^\top \phi)(x) \approx \mathcal{A}_k(\lambda^\top \phi(x-\cdot) + \|\cdot\|_2^2/2)$ with $\mathcal{A}_k$ the $k$\textsuperscript{th} iteration of an algorithm $\mathcal{A}$ initialized at $0$. Ultimately, we can also consider approximating the Moreau gradient of the whole dynamical system as in~\eqref{eq:chain} by an optimization subroutine as presented in Sec.~\ref{sec:moreau} for generic functions. The decomposition of the approximation of the Moreau gradients at different scales provides then different oracles. 

Several subroutines can be used to approximate the Moreau gradients. 
As shown in Lemma~\ref{lem:stepsize}, under appropriate smoothness assumptions, we can select the scaling parameters of the intermediate Moreau gradient computations of a composition in such a way that the subproblem is strongly convex. In that case, we can use any first-order order method to approximate the Moreau gradient. For example, a gradient descent on the subproblem can approximate the Moreau gradient up to $\varepsilon$ accuracy in at most $O(\log\varepsilon)$ iterations. In the extreme case, if one considers a single gradient step $\GD^{(1)}$ from $y=0$ to approximate the Moreau gradient or its augmented version, we get
\begin{align*}
	\widehat \nabla \env(f)(z) & = \GD^{(1)}(\lambda^\top f(z-\cdot) +\|\cdot\|_2^2/2)  = \nabla f(z)\lambda \\
	\widehat \nabla \anv_\reglag(f)(z;\lambda)  & = \GD^{(1)}(\reglag \|f(z - \cdot) - f(z) + \lambda/\reglag\|_2^2 + \|\cdot\|_2^2) = \nabla f(z)\lambda,
\end{align*}
such that the classical gradient vector product rule can be seen as an approximation of the chain rules presented for the Moreau gradient and its augmented version. By varying the procedures used to approximate the Moreau gradients or the augmented Moreau gradients of the intermediate computations, we can recast  back-propagation schemes proposed earlier by~\citet{frerix2018proximal} and~\citet{lee2015difference} as explained in Sec.~\ref{sec:rel_algos}.

\subsection{Computational Complexities}
In terms of memory usage, our approximations of the Moreau gradients have the same cost as the ones of gradient back-propagation as implemented in modern machine learning frameworks such as  Pytorch~\citep{paszke2017automatic} and~\citep{tensorflow2015-whitepaper} which also proceed by keeping in memory the intermediate inputs and the code implementing the function. We detail here the computational cost of an oracle based on Moreau gradients as in Algo.~\ref{algo:ambp_backward}, the computational complexity of an oracle based on augmented Moreau gradients is analogous when  using subroutines to approximate the solutions of the subproblems. 

In terms of computational complexity, our approach differs depending on the number of iterations used by the optimization subroutine. To state the computational complexities, we denote by $\mathcal{T}(f)$, $\mathcal{T}(\nabla f)$ the computational complexity of evaluating a function $f$ and its gradient respectively, and we denote by $\mathcal{T}(\nabla f(x))$ the computational complexity of computing a product of the form $\nabla f(x) \lambda$, i.e., we identify $\nabla f(x)$ to the corresponding linear function.
\begin{enumerate}[nosep]
	\item For a sequence of compositions of the form $f=  \varphi_\horizon \circ \ldots\varphi_1: \reals^{d_0} \rightarrow \reals^{d\horizon}$, given $x \in \reals^{d_0}$ and denoting $x_t = \phi_t(x_{t-1})$ for $t \in \{1, \ldots, \horizon\}$ with $x_0 = x$,
	\begin{enumerate}[nosep]
		\item the computational cost of computing $\nabla (\lambda ^\top f)(w)$ by gradient back-propagation is of the order of
		\[
		\mathcal{T}_{\BP} = \sum_{t=1}^\horizon \left[\mathcal{T}(\varphi_t)  + \mathcal{T}(\nabla \varphi_t) + \mathcal{T}(\nabla \varphi_t(x_{t-1}))\right],
		\]
		\item the computational complexity of approximating the Moreau gradient $\nabla \env (\lambda^\top f)$ by chain rule applications as in~\eqref{eq:mbp_backprop_detailed} with a first order optimization subroutine $\mathcal{A}_t$ to approximate $\nabla \env(\lambda_t^\top \varphi_t)(x_t)$ as $\mathcal{A}(\lambda_t^\top \varphi_t(x_t - \cdot) + \|\cdot\|_2^2/2)$ is of the order of 
		\[
		\mathcal{T}_{\MBP} = \sum_{t=1}^\horizon \left[\mathcal{T}(\varphi_t)  + K_{\mathcal{A}_t}\left(\mathcal{T}(\nabla \varphi_t) + \mathcal{T}(\nabla \varphi_t(x_{t-1}))\right)\right],
		\]
		where $K_{\mathcal{A}_t}$ is the number of oracle calls of the algorithm $\mathcal{A}_t$ initialized at 0 to output  $\mathcal{A}(\lambda_t^\top \varphi_t(x_t - \cdot) + \|\cdot\|_2^2/2)$. If $\mathcal{A}_t$ is a gradient descent $\GD$, $K_{\mathcal{A}_t}$ is the number of iterations of $\GD$. 
	\end{enumerate}	
	\item For a dynamical system as in~\eqref{eq:chain}, i.e., a function $f:\reals^p \rightarrow\reals^{d_\horizon}$ such that  for $w=(w_1, \ldots, w_{\horizon})$, we have $f(w) = x_\horizon$ with $x_t = \dyn_t(w_t, x_{t-1})$ for $t \in \{1, \ldots, \horizon\}$ with $x_0 $ fixed,
	\begin{enumerate}
		\item the computational complexity of computing $\nabla (\lambda^\top f)(w)$ by gradient back-propagation is of the order of 
		\[
		\mathcal{T}_{\BP} = \sum_{t=1}^\horizon \left[\mathcal{T}(\dyn_t) + \mathcal{T}(\nabla \dyn_t) + \mathcal{T}(\nabla \dyn_t(w_t, x_{t-1})) \right],
		\]
		\item the computational complexity of approximating the Moreau gradients $(\nabla \env(\lambda^\top f_{w, t})(w_t))_{t=1}^\horizon$ by Algo.~\ref{algo:mbp_backward}, for $f_{w, t}$ defined in~\eqref{eq:shift_coord}, with first-order optimization subroutines $\mathcal{A}_{t}, \mathcal{B}_{t}$ to approximate $\nabla \env(\lambda_t^\top \phi_t(w_t, \cdot))(x_{t-1})$ and $\nabla \env(\lambda_t^\top \phi_t(\cdot, x_{t-1}))(w_t)$ is of the order of 
		\begin{align*}
	\hspace{-20pt}	
	\mathcal{T}_{\MBP} = \sum_{t=1}^\horizon \big[\mathcal{T}(\dyn_t) 
	& + K_{\mathcal{A}_t}\left(\mathcal{T}(\nabla_{x_{t-1}} \dyn_t) + \mathcal{T}(\nabla_{x_{t-1}} \dyn_t(w_t, x_{t-1}))\right) \\
	& + K_{\mathcal{B}_t} \left(\mathcal{T}(\nabla_{w_t} \dyn_t) + \mathcal{T}(\nabla_{w_t} \dyn_t(w_t, x_{t-1})) \right) \big],
		\end{align*}
	with $K_{\mathcal{A}_t}, K_{\mathcal{B}_t}$ defined as above and $\nabla_{x_{t-1}} \dyn_t, \nabla_{w_t} \dyn_t$ denoting partial gradient evaluations. 
	\end{enumerate}
\end{enumerate}
In practice, we use a quasi-Newton algorithm to approximate the Moreau gradients that exhibit fast convergence with few oracles evaluations for each sub-problem.

\section{Related Algorithms}\label{sec:rel_algos}
We present here two earlier back-propagation schemes namely target propagation~\citep{lee2015difference} and proximal back-propagation~\citep{frerix2018proximal} that can be cast as approximate computations of the Moreau envelope. The root of the different back-propagation schemes can be found in the formulation of the Moreau gradients in terms of the penalized formulation of the problem. 

\subsection{Target Propagation}
We start by recalling the rationale of target propagation as presented by, e.g.,~\citet{lee2015difference}. Earlier implementations of target propagation schemes have been developed by~\citet{lecun1986learning, rohwer1990moving, bengio2013estimating} and more recently by~\citet{manchev2020target, fairbank2022deep}. 

The  idea of target propagation is to compute virtual targets $\target_t$ for each layer $\varphi_{t} = \phi_t(w_t, \cdot)$ $t=1, \ldots, \horizon$ such that if the layers were able to match their corresponding target, i.e., 
$ \varphi_t(x_{t-1}) \approx \target_t$, 
the  objective  would decrease. The final target $\target_\horizon$ is computed  as a  gradient step on the loss, i.e., the outer function $h$ w.r.t. $x_\horizon$.  The targets are then back-propagated using an approximate inverse
$\varphi_{t}^{\invsign}$ of each layer $\varphi_t$ and the parameters $w_t$ are updated by a gradient step on $\| \phi_t(\cdot, x_{t-1}) - \target_t\|_2^2$.  The initial rationale of target propagation was to propagate the targets through the approximate inverses directly, i.e., using $z_{t-1} = \varphi_t^\invsign(z_t)$. However, this approach has not been successful empirically and was modified by introducing the\emph{ difference target propagation formula}  $z_{t-1} = x_{t-1} + \varphi_t^\invsign(z_t)  - \varphi_t^\invsign(x_t)$ for $x_t$ the intermediate states of the layers, see~\citet{lee2015difference}.
 
Formally, consider a feed-forward network $\phi:\reals^{d_0} \times \reals^p \rightarrow \reals^{d_\horizon}$ such that for $x_0 \in \reals^{d_0}, w = (w_1, \ldots, w_\horizon) \in \reals^p$, the output of the network on $x_0$ is $\phi(w, x_0) = x_\horizon$ with $x_\horizon$ given by passing $x_0$ through $\horizon$ layers $\phi_t$ such that $x_t = \phi_t(w_t, x_{t-1})$ for $t=1, \ldots, \horizon$. Given a set of parameters $w$ and an input $x_0$, target propagation with a difference target propagation formula proceeds by computing $x_t$ for $t=1, \ldots, \horizon$ in a forward pass, then by performing a backward pass that propagates targets as
\begin{align}\label{eq:tp_backprop}
	z_\horizon = x_\horizon - \auxstepsize_\horizon\nabla h(x_\horizon), \quad
	z_{t-1} = x_{t-1} + \varphi_t^\invsign(z_t)  - \varphi_t^\invsign(x_t) \ \mbox{for} \ t \in \{2, \ldots, \horizon\}, 
\end{align}
and by updating the parameters as
\begin{align}\label{eq:tp_update}
	w_t  \leftarrow w_t - \stepsize_t \nabla (\| \phi_t(\cdot, x_{t-1}) - \target_t\|_2^2)(w_t) \ \mbox{for} \ t\in \{1, \ldots, \horizon\}.
\end{align}
The difference target propagation formula can be interpreted as an approximation of the linearization of the approximate inverse, which itself can be interpreted as using the inverse of the gradient of the layer, namely $z_{t-1} - x_{t-1} = \varphi_t^\invsign(z_t)  - \varphi_t^\invsign(x_t) \approx \nabla \varphi_t^\invsign (x_t)^\top (z_t -x_t) \approx( \nabla \varphi_t(x_{t-1}))^{-1}(z_t -x_t)$ using that $\varphi_t^\invsign \circ \varphi_t(x_{t-1}) \approx x_{t-1}$ such that $\nabla \varphi_t^\invsign(x_t) \approx ( \nabla \varphi_t(x_{t-1}))^{-1}$. The back-propagation scheme presented in Eq.~\eqref{eq:tp_backprop} and Eq.~\eqref{eq:tp_update} can then be rewritten in terms of the displacements $\delta_t = z_t - x_t$ as 
\begin{align*}\label{eq:tp_reformulated}
	\delta_\horizon & =  \auxstepsize_\horizon\nabla h(x_\horizon), \quad 
	\delta_{t-1} \approx( \nabla \varphi_t(x_t))^{-1} \delta_t \hspace{20pt}
	\mbox{for} \ t \in \{\horizon, \ldots, 2\}\\
	w_t & \leftarrow  w_t - \stepsize_t g_t, \hspace{19pt} 
	g_t = \nabla_{w_t} \phi_t(w_t, x_{t-1}) \delta_t \ 
	\mbox{for} \ t \in \{\horizon, \ldots, 1\}.
\end{align*}

Consider now again the back-propagation scheme based on augmented Moreau gradients as presented in Algo.~\ref{algo:ambp_backward}.  If the Moreau gradient on the outer function $h$ is approximated with a gradient step, we retrieve $\lambda_\horizon = \auxstepsize_\horizon  \nabla h(x_\horizon)$. The propagation of the variables $\lambda_t$ amounts to computing
\[
\lambda_{t-1} \approx \argmin_{y_{t-1}}  \|\phi_t(w_t, x_{t-1} - y_{t-1}) - \phi_t(w_t, x_{t-1}) + \lambda_t\|_2^2 +  \rho\|y_{t-1}\|_2^2,
\]
for $\rho= (\reglag \auxstepsize_{t-1})^{-1}$ some regularization parameter. By approximating the solution of the above problem using a Gauss-Newton step, i.e., by approximating $\phi_t(w_t, x_{t-1} - y_{t-1}) = \phi_t(w_t, x_{t-1})  - \nabla_{x_{t-1}} \phi_t(w_t, x_{t-1}) y_{t-1}$ and solving the resulting problem, we get $\lambda_{t-1} = (\nabla \varphi_t(x_t) \nabla \varphi_t(x_t)^\top + \rho\idm)^{-1} \nabla \varphi_t (x_t)  \lambda_t$ and for $\rho\ll 1$, we have $\lambda_{t-1} \approx ( \nabla \varphi_t(x_t))^{-1}\lambda_t$. In other words, the back-propagation with augmented Moreau gradients uses regularized inverses of the gradients rather than using the inverses of the gradients directly. Finally, if we consider using a gradient step to approximate the computation of $g_t$ in Algo.~\ref{algo:ambp_backward}, we get $g_t = \stepsize_t \reglag \nabla_{w_t} \phi_t(w_t, x_{t-1})$. Hence target propagation with a difference target propagation formula can be seen as a possible implementation of a backward pass using augmented Moreau gradients. Compared to target propagation, Algo.~\ref{algo:ambp_backward} introduces an additional regularization that may stabilize the back-propagation scheme. This regularization can be interpreted as a stabilization procedure akin to previous heuristics used to implement target propagation~\citep{roulet2021target}.

\subsection{Proximal Back-propagation}
\citet{frerix2018proximal} considered feed-forward networks as defined above for target propagation. They further decompose the layers as $\phi_t(w_t, x_{t-1}) = a_t(b_t(w_t, x_{t-1}))$, where $a_t$ is a non-linear activation function applied element-wise and $b_t$ is a bilinear function such as a matrix-vector product with $w_t$ containing the matrix and $x_{t-1}$  being the vector. After noticing that a gradient step using gradient back-propagation can be seen as a block coordinate gradient pass on the penalized formulation of the problem, \citet{frerix2018proximal} propose to blend gradient back-propagation and a proximal step on the parameters of the network. After reparameterizing the approach in terms of the displacements incurred by the block coordinate gradient pass, their algorithm can be written as
\begin{align*}
	\lambda_\horizon & = \nabla h(x_\horizon), \\
	\mu_{t} & = \nabla a_t(b_t(w_t, x_{t-1}) \lambda_t, \quad \lambda_{t-1} = \nabla_{x_{t-1}} b_t(w_t, x_{t-1}) \mu_{t}  \quad \mbox{for} \ t \in \{1, \ldots, \horizon\}\\
	w_t & \leftarrow w_t - \widehat \nabla \anv_{\reglag_t}(b_t(\cdot, x_{t-1}))(w_t; \reglag_t\mu_t),
\end{align*}
where $\reglag_t$ is some stepsize parameter and $\widehat \nabla \anv_{\reglag_t}(b_t(\cdot, x_{t-1}))(w_t; \reglag_t\mu_t)$ is computed by a conjugate gradient method. 

We retrieve here an implementation of the backward pass using augmented Moreau gradients where the latter are approximated by a single gradient step for all intermediate computations except for the computation of the bilinear operation where a conjugate gradient method is used. 
The proximal back-propagation scheme of~\citet{frerix2018proximal} is then an instance of the schemes presented in Algo.~\ref{algo:gen_forward} and Algo.~\ref{algo:gen_backward}  where we compose the usual $\GBP$ rule for some operations and one of the proposed $\MBP$ or $\IBP$ rules for some other operations. 

\section{Optimization with Moreau Gradients}\label{sec:oracle}
Algos.~\ref{algo:gen_forward} and~\ref{algo:gen_backward} summarize a generic scheme to compute oracles on an objective. These oracles can be adapted to the form of the objective, namely, finite-sum objectives by considering mini-batches and can be used within different optimization algorithms as shown in this section. 

\subsection{Stochastic Setting}
For deep learning problems, the objective blends a dynamical structure of the form~\eqref{eq:chain} and a finite sum. Namely, the training objective is of the form 
\begin{equation}\label{eq:learning_pb}
\min_{w \in \reals^p} \sum_{i=1}^{n} \mathcal{L}(y^{(i)}, \psi(w, x^{(i)})),
\end{equation}
where $(x^{(i)}, y^{(i)})$ are samples of input-output pairs, $\mathcal{L}(y, \hat y)$ is the loss incurred by predicting label $\hat y$ instead of the true label $ y$ and $\psi(w, x^{(i)}) = \psi^{(i)}(w)$ is a deep network composed of $\horizon$ layers $\phi_t$ with parameters $w_t$ for $t={1, \ldots, \horizon}$ encompassed in $w=(w_1, \ldots, w_\horizon)$  such that
\begin{align}
 \psi^{(i)}(w) & = x_{\horizon}  \nonumber\\
 x_t & = \phi_t(w_t, x_t), \quad  \mbox{for} \ t \in \{1, \ldots \horizon\} \ x_0 = x^{(i)}. \label{eq:deep_dyn}
\end{align}
For a multi-layer perceptron with ReLU activation functions taking vectors $x^{(i)}$ as inputs, the layers can be written as $\phi_t(w_t, x_{t-1}) = r(W_tx_{t-1} + b_{t-1})$ with $r$ the element-wise application function of the ReLU activation $\relu(z) = \max\{z, 0\}$ and $w_t = (W_t, b_t)$ composed of the weight matrix and the offset parameters $b_t$. 

The objective in~\eqref{eq:learning_pb} can be written as a composition of the form $h\circ f$ by defining 
\[
h(z) = \frac{1}{n} \sum_{i=1}^n \mathcal{L}(y^{(i)}, z^{(i)}) \quad f(w) = (\psi(w, x^{(1)}), \ldots, \psi(w, x^{(n)})),
\]
such that $f$ keeps a dynamical structure as in~\eqref{eq:deep_dyn} by simply concatenating the dynamical structures of the functions $\psi^{(i)}$ defined previously. However, while $\psi^{(i)}$ mapped the parameters to, e.g., a real number, the function $f$ maps the parameters to a vector of size $n$. To alleviate the potential increased complexity done by considering concatenated dynamics, we consider computing Moreau gradients on mini-batches, i.e., functions of the form
\begin{equation}\nonumber
	\obj_S(\chain_S(\ctrls)) = \frac{1}{|S|}\sum_{i\in S} \loss(\labl^{(i)}, \chainnet(\ctrls, x^{(i)})),
\end{equation}
where  $S$ is a mini-batch, $h_S$ is defined analogously as above and $\chain_S$ is the concatenation of the functions $\ctrls\mapsto\chainnet(\ctrls, x^{(i)})$, i.e., $f_S$ keeps a dynamical structure but maps now to a vector of  $|S|$ predictions. 

\subsection{Composite Objectives}
In several applications, the objective is a sum of an objective possessing a dynamical structure and a regularization term, i.e, an objective of the form $\outfunc(\infunc(w)) + \reg(w)$, where $\reg(w)$ decomposes along the parameters of the dynamical system defined by $f$, such that $\reg(w) = \sum_{t=1}^\horizon \reg(w_i)$ for $w=(w_1, \ldots, w_\horizon)$, and, e.g., $\reg(w) = \|w\|_2^2$. In that case, we consider computing an approximation of the Moreau gradient of the sum by keeping the forward and backward passes presented in Algo.~\ref{algo:gen_forward} and Algo.~\ref{algo:gen_backward} and simply change the computation of the oracle by introducing the regularizer in the minimization. For the approximation of the Moreau gradients, the oracle direction is then given as 
\[
g_t \approx \argmin_{v_t} \stepsize_t\adjvar_t^\top \dyn_t(w_t-v_t, x_{t-1}) + \frac{1}{2}\|v_t\|_2^2 + \reg(w_t-v_t),
\]  
instead of $g_t \approx \argmin_{v_t}  \stepsize_t\adjvar_t^\top \dyn_t(w_t-v_t, x_{t-1}) + \frac{1}{2}\|v_t\|_2^2$ on line~\ref{line:oracle_mbp} in Algo.~\ref{algo:mbp_backward}, where the approximation is done with an optimization algorithm.
\subsection{Optimization Algorithms}
Given the outputs returned by procedures of the form Algo~\ref{algo:gen_forward} and~\ref{algo:gen_backward}, we can consider various update rules for the parameters. For example, one can consider using the output directions $g_t$ of Algo~\ref{algo:gen_backward} to update the parameters $w_t$ as $w_t \leftarrow w_t -g_t$ as explained in Sec.~\ref{sec:moreau}. Alternatively, one can plug the directions $g_t$ in an optimization method used for stochastic optimization such as Adam~\citep{kingma2015adam} or a Stochastic Gradient Descent with momentum. For example, by plugging the directions in SGD with Nesterov momentum without dampening as implemented in PyTorch~\citep{paszke2017automatic}, the updates take then the form $w_t \leftarrow w_t -  \gamma (g_t + \nu b_t)$ where $b_t =  \sum_{s=1}^t \nu^{t-s} g_s$ is computed and updated along the iterations. The usual implementation in PyTorch takes $g_t$ to be stochastic estimates of the gradients. We can consider using stochastic estimates of the Moreau gradients instead.

\paragraph{Acknowledgments}
This work was supported by NSF DMS-2023166, NSF CCF-2019844, NSF DMS-1839371, CIFAR-LMB, and faculty research awards.

\bibliography{soft_bp}
\bibliographystyle{agsm}

\clearpage
\appendix

\section*{\LARGE Appendix}
The Appendix is organized as follows.
\begin{enumerate}[nosep]
	\item Appendix~\ref{app:moreau} presents convergence results with approximate Moreau gradients.
	\item Appendix~\ref{app:chain_rule} details the derivations of chain rules for Moreau gradients.
	\item Appendix~\ref{app:helper} presents helper lemmas for both convergence results and derivations of the chain rules.
	\item Appendix~\ref{app:exp} presents experimental details on the numerical illustrations.
\end{enumerate}

\section{Convergence Guarantees with Moreau Gradients}\label{app:moreau}
\subsection{Optimization with Moreau Gradients}
In this section, we recall convergence guarantees associated with an approximate Moreau gradient descent and its augmented version and compare them to the analysis of gradient descent.

\paragraph{Gradient descent}
Consider minimizing a function $f:\reals^d \rightarrow \reals$. A standard gradient descent with stepsize $\stepsize>0$ reads
\begin{equation}\label{eq:gd_app}
w^{(k+1)} = w^{(k)} - \stepsize \nabla f(w^{(k)}).
\end{equation}
The maximal stepsize  that ensures a decrease of the objective values is a priori bounded by the inverse of the Lipschitz continuity parameter, $\smooth_f$,  of the gradients of the objective. By taking $\stepsize \leq  1/\smooth_f$, we can guarantee that, after $K$ iterations, we get a point that is nearly stationary, i.e, we have that
\[
\min_{k\in\{0, \ldots, K-1\}}  \|\nabla f(w^{(k)})\|_2^2 \leq \frac{2(f(w^{(0)})- f^*)}{\stepsize K},
\]  
where $f^* = \min_{w \in \reals^d} f(w)$, see, e.g.,~\citep{nesterov2013introductory}.

\paragraph{Approximate Moreau gradient descent}
The maximal stepsize to define the Moreau gradient is the maximal value of $\bar \stepsize$ such that for any $w\in \reals^d$, $v \mapsto  \bar \stepsize f(w-v) + \|v\|_2^2/2$ is convex. If, for example, $f$ is convex, then $\bar \stepsize = +\infty$. 
On the other hand, if $f$ is $\smooth_f$-smooth, that is, with $\smooth_f$-Lipschitz continuous gradients, then $\bar \stepsize \geq 1/\smooth_f$ as recalled in Lemma~\ref{lem:smooth_strg_cvx}.
 As explained in Sec.~\ref{sec:moreau}, the Moreau gradient may not be available in closed form but it may be approximated up to any accuracy by some subroutine. We consider then an approximate Moreau gradient descent with stepsize $\stepsize$ as the sequence of iterates
\begin{equation}\label{eq:mgd_app}
w^{(k+1)} = w^{(k)} - \widehat \nabla \env(\stepsize f) (w^{(k)}), \quad \mbox{where}\  \|\widehat \nabla \env(\stepsize f)(w^{(k)}) - \nabla \env(\stepsize f)(w^{(k)})\|_2 \leq \stepsize \varepsilon_k.
\end{equation}
For $\stepsize \leq \bar \stepsize/2$, the condition in~\eqref{eq:mgd_app} can be verified by finding $ \widehat \nabla \env(\stepsize f) (w)  = v$ such that $\|\stepsize \nabla f(w-v) + v\| \leq \stepsize\varepsilon_k$ since for a 1-strongly function such as $g:v \rightarrow \stepsize f(w-v) + \|v\|_2^2/2$, we have $\|v^*- v\|_2 \leq \|\nabla g(v)\|_2$ for $v^*=\argmin_v g(v) =  \nabla \env(\stepsize f)(w)$~\citep[Theorem 2.1.10]{nesterov2013introductory}.
The potentially larger stepsize taken by the Moreau gradient is balanced by the inaccuracy of the oracle as recalled in the following lemma.

\begin{lemma}\label{lem:mbp_cvg}
	Consider $f:\reals^d \rightarrow \reals$ and $\bar \stepsize>0$ such that for any $w\in \reals^d$, $v \mapsto  \bar \stepsize f(w-v) + \|v\|_2^2/2$ is convex. The iterates of an approximate Moreau gradient descent~\eqref{eq:mgd_app} with stepsize $\stepsize \leq \bar \stepsize/2$ satisfy 
	\begin{align*}
		\min_{k\in\{0, \ldots, K-1\}} \frac{\|\nabla \env(\stepsize f)(w^{(k)})\|_2^2}{\stepsize^2} \leq  {\frac{2(f(w^{(0)})  -f^*)}{\stepsize K}}  +  \frac{1}{K}\sum_{k=0}^{K-1} \varepsilon_k^2,
	\end{align*}
	where $f^* = \min_{w \in \reals^d} f(w)$.
\end{lemma}
\begin{proof}
 Consider the function $f_\stepsize: w \mapsto \env(\stepsize f)(w)/\stepsize$. We have that $\nabla f_\stepsize(w) = \nabla \env(\stepsize f)(w)/\stepsize$, hence the iterates~\eqref{eq:mgd_app} are an approximate gradient descent on $f_\stepsize$ with stepsize $\stepsize$, namely, they can be written $w^{(k+1)} = w^{(k)} - \stepsize \widehat \nabla f_\stepsize(w^{(k)})$ with $ \|\widehat \nabla f_\stepsize(w^{(k)})- \nabla f_\stepsize(w^{(k)})\|_2\leq \varepsilon_k$. 
 Moreover, recall that for $\stepsize \leq \bar\stepsize/2$, the Moreau gradients are $1$-Lipschitz continuous, hence $f_\stepsize$ is $1/\stepsize$-smooth. Hence, by Lemma~\ref{lem:approx_grad_cvg}, we have that after $K$ iterations, 
 \begin{equation}\nonumber
\min_{k\in\{0, \ldots, K-1\}}\|\nabla f_\stepsize(w^{(k)})\|_2^2 \leq  {\frac{2(f_\stepsize(w^{(0)})  - f_\stepsize^*)}{\stepsize K}}  +  \frac{1}{K}\sum_{k=0}^{K-1} \varepsilon_k^2,
 \end{equation}
 where $f_\stepsize^*= \min_{w \in \reals^d} f_\stepsize(w))$.
Now, by definition of $f_\stepsize$, we have that $f_\stepsize(w) \leq f(w)$ for all $w \in \reals^d$ and $f_\stepsize(w) \geq  \min_{w' \in \reals^d} f(w')$ for all $w \in \reals^d$ which concludes the claim. 
\end{proof}

In Lemma~\ref{lem:mbp_cvg}, we use $\|\nabla \env(\stepsize f)(w)/\stepsize\|_2$ as a measure of stationarity of the problem. Indeed, the nullity of $\nabla \env(\stepsize f)(w)$ is a first-order necessary optimality condition for the problem, since if $w^* \in \argmin_{w\in \reals^d} f(w)$ then $\nabla \env(\stepsize f)(w^*) = 0$. As $\stepsize$ is used as a stepsize for the method we consider the normalized quantity $\|\nabla \env(\stepsize f)(w)/\stepsize\|_2$. As shown in the following lemma, the quantity $\|\nabla \env(\stepsize f)(w)/\stepsize\|_2$ can be translated as a criterion to near stationarity of a point in the classical sense, i.e., in terms of the norm of the gradient. 

\begin{lemma}\label{lem:approx_stat_mor}
	Consider the assumptions of Lemma~\ref{lem:mbp_cvg}. If a point $w \in \reals^d$ satisfies 
	$
	\|\nabla \env(\stepsize f)(w)/\stepsize\|_2 \leq \varepsilon,
	$
	for any $\stepsize < \bar \stepsize$,
	then $w$ is close to a point $w^*$ that is nearly stationary, in the sense that 
	\begin{align*}
		\|w - w^*\|_2 \leq \stepsize \varepsilon, \qquad \|\nabla f(w^*)\|_2 \leq \varepsilon.
	\end{align*}
\end{lemma}
\begin{proof}
	Consider $w^* = w-\nabla\env(\stepsize f)(w)$. By  definition $\|w - w^*\|_2 = \|\nabla \env(\stepsize f)(w)\|\leq \stepsize\varepsilon$ and $w^*= \argmin_{w' \in \reals^d} f(w') + \|w-w'\|_2^2/(2\stepsize)$  such that $\nabla f(w^*) =( w-w^*)/\stepsize =  \nabla \env(\stepsize f)(w)/\stepsize$ which gives the claim. 
\end{proof}

\paragraph{Augmented approximate Moreau gradient descent}
Recall that we defined the augmented Moreau gradient for any $\stepsize$ and $\reglag\geq 0$ as 
\begin{align*}
	\nabla \anv_\augparam( f)(w; \stepsize) & = \argmin_{v \in \reals^\dimin} \stepsize f(w{-}v) {+} \frac{1}{2} \|v\|_2^2 {+} \frac{\augparam}{2} (f(w{-}v){-}f(w))^2
	= \argmin_{v \in \reals^\dimin} \augparam (f(w{-}v) {-} f(w) {+} \stepsize/\augparam)^2 {+}  \|v\|_2^2.
\end{align*}
As mentioned in Sec.~\ref{sec:moreau}, the augmented Moreau gradient is defined for any $\stepsize$ provided that $\reglag$ is positive. However, the augmented Moreau gradient is not necessarily Lipschitz-continuous for any $\stepsize, \reglag>0$ which prevents direct use of Lemma~\ref{lem:approx_grad_cvg}. We present simple sufficient conditions to ensure that the augmented Moreau gradient is Lipschitz continuous.
Consider $\bar \reglag$ such that for any $w\in \reals^\dimin$, $v \mapsto \bar \reglag (f(w-v) -f(w))^2/2 + \|v\|_2^2/2$ is convex and $\bar \stepsize$ such that $v \mapsto \bar \stepsize f(w-v) + \|v\|_2^2/2$ is convex. For any $\reglag \leq \bar \reglag/2$ and $\stepsize \leq \bar \stepsize /2$, the augmented Moreau gradient is uniquely defined and $\max\{\stepsize/(\bar \stepsize/2- \stepsize), 1\}$-Lipschitz-continuous, that is, $1$-Lipschitz continuous for $\stepsize \leq \bar \stepsize/4$. For $f$ that is $\smooth_f$-smooth, we already know that $\bar \stepsize \geq \smooth_f$. On the other hand, a lower bound for $\bar \reglag$ is  $1/(2 \smooth_f b_f)$ for a $\smooth_f$-smooth and $b_f$-bounded function  $f$. Equipped with appropriate parameters, we can then analyze an approximate augmented Moreau gradient descent of the form 
\begin{equation}\label{eq:amp_app}
w^{(k+1)} = w^{(k)} - \widehat \nabla \anv_\reglag(f)(w^{(k)}; \stepsize), \qquad \mbox{where} \ \| \widehat \nabla \anv_\reglag(f)(w^{(k)}; \stepsize) -   \nabla \anv_\reglag(f)(w^{(k)}; \stepsize)\|_2 \leq \stepsize \varepsilon_k,
\end{equation}
as shown in the following lemma. As for the Moreau gradient, condition~\eqref{eq:amp_app} can be verified from the gradient of $g: v \mapsto \stepsize f(w-v) + \|v\|_2^2/2 + \reglag (f(w-v) - f(w))^2/2$ provided that $g$ is 1-strongly convex which is the case for  $\reglag \leq \bar \reglag/2$ and $\stepsize \leq \bar \stepsize /4$.
\begin{lemma}\label{lem:approx_ampb}
	Consider $f:\reals^d \rightarrow \reals$, $\bar \reglag$, $\bar \stepsize>0$ such that for any $w\in \reals^d$, $v \mapsto  \bar \stepsize f(w-v) + \|v\|_2^2/2$ and $v \mapsto \bar \reglag (f(w-v) -f(w))^2/2 + \|v\|_2^2/2$ are convex. The iterates of an approximate augmented Moreau gradient descent~\eqref{eq:mgd_app} with parameter $\reglag \leq \bar\reglag/2$ and stepsize $\stepsize \leq \bar \stepsize/4$ satisfy
\begin{align*}
	\min_{k\in\{0, \ldots, K-1\}} \frac{\| \nabla \anv_\reglag(f)(w^{(k)}; \stepsize) \|_2^2}{\stepsize^2} \leq  {\frac{2(f(w^{(0)})  -f^*)}{\stepsize K}}  +  \frac{1}{K}\sum_{k=0}^{K-1} \varepsilon_k^2,
\end{align*}
where $f^* = \min_{w \in \reals^d} f(w)$.
\end{lemma}
\begin{proof}
	Consider the function $f_{\stepsize, \reglag}: w\rightarrow  \anv_\reglag(f)(w; \stepsize)/\stepsize$. We have that the iterates~\eqref{eq:amp_app} can be written as $w^{(k+1)} = w^{(k)} - \stepsize \widehat \nabla f_{\stepsize, \reglag}(w^{(k)})$ with $\|\widehat \nabla f_{\stepsize, \reglag}(w) -  \nabla f_{\stepsize, \reglag}(w)\|_2 \leq \varepsilon_k$. Given the assumptions, we have that $\nabla f_{\stepsize, \reglag}$ is $1/\stepsize$-Lipschitz-continuous. Hence, by Lemma~\ref{lem:approx_grad_cvg}, we have that 
	 \begin{equation}\nonumber
		\min_{k\in\{0, \ldots, K-1\}}\|\nabla f_{\stepsize, \reglag}(w^{(k)})\|_2^2 \leq  {\frac{2(f_{\stepsize, \reglag}(w^{(0)})  - f_{\stepsize, \reglag}^*)}{\stepsize K}}  +  \frac{1}{K}\sum_{k=0}^{K-1} \varepsilon_k^2,
	\end{equation}
 where $f_{\stepsize, \reglag}^*= \min_{w \in \reals^d} f_{\stepsize, \reglag}(w)$.
	Moreover, $f_{\stepsize, \reglag}(w) \leq f(w)$ and $f_{\stepsize, \reglag}(w) \geq \min_{w'\in \reals^\dimin}f(w')$, which concludes the proof. 
\end{proof}
As for the Moreau gradient descent, we used the scaled augmented Moreau gradient as a measure of stationarity as it defines a necessary optimality condition. On the other hand, Lemma~\ref{lem:approx_stat_mor} can be adapted to this case.

\begin{lemma}
	Consider the assumptions of Lemma~\ref{lem:approx_ampb} and assume in addition that $f$ is $\ell_f$-Lipschitz continuous. For any $\stepsize \leq \bar \stepsize/2$ and $\reglag \leq \bar \reglag/2$, if a point $w\in \reals^\dimin$ satisfies $\|\nabla \anv_\reglag(f)(w; \stepsize)/\stepsize\|_2\leq \varepsilon$ with $\varepsilon \leq 1/(2\reglag\lip_f)$, then $w$ is close to a point $w^*$ that is nearly stationary, in the sense that 
	\[
	\|w-w^*\|_2 \leq \stepsize\varepsilon, \qquad  \|\nabla f(w^*) \|_2 \leq 2 \varepsilon.
	\]
\end{lemma}
\begin{proof}
Consider $w^* = w-\nabla\anv_\reglag(f)(w; \stepsize)$. By  definition $\|w - w^*\|_2 = \|\nabla \anv_\reglag(f)(w; \stepsize)\|\leq \stepsize\varepsilon$ and $w^*= \argmin_{w' \in \reals^d} f(w') + \|w-w'\|_2^2/(2\stepsize) + \reglag (f(w') -f(w))^2/(2\stepsize)$  such that $\nabla f(w^*) = ( w-w^*)/\stepsize + \reglag \nabla f(w^*) (f(w) - f(w^*))/\stepsize$ and 
$\|\nabla  f(w^*) \|_ 2 \leq \|w-w^*\|_2/\stepsize + \reglag \lip_f\|\nabla f(w^*)\|_2 \|w-w^*\|_2/\stepsize$. 
Rearranging the terms, using that $\|w-w^*\|_2/\stepsize \leq \varepsilon$ and $\varepsilon \leq 1/(2\reglag\lip_f)$ gives the claim. 
\end{proof}

\paragraph{Summary}
The following corollary summarizes the results presented in this section. The complexity bounds of gradient descent and an approximate (augmented) Moreau gradient descent are similar in the sense that the number of iterations to reach a point that is $\varepsilon$ stationary in terms of gradient norm is of the order of $\varepsilon^2$ for both methods. However, the constants may differ since the stepsize taken by a Moreau gradient descent is a priori larger than the stepsize taken for a standard gradient descent. 

\begin{corollary}
	Consider $f:\reals^d \rightarrow \reals$. 
	The number of iterations of a gradient descent~\eqref{eq:gd_app} to reach a point $\hat w \in \reals^d$ that is $\varepsilon$ stationary, i.e., $\|\nabla f(\hat w)\|_2 \leq \varepsilon$ is at most
	\[
	K \leq \frac{2 (f(w^{(0)})- f^*)}{\stepsize \varepsilon^2},
	\]
	provided that the stepsize satisfies $\stepsize \leq  1/\smooth_f$, with $\smooth_f$ the smoothness parameter of $f$.
	
	The number of iterations of an approximate Moreau gradient descent~\eqref{eq:mgd_app} with constant approximation error $\varepsilon/\sqrt 2$ to output a point that is $\stepsize\varepsilon$ close to a $\varepsilon$ stationary point is at most
	\[
	K \leq \frac{4 (f(w^{(0)})- f^*)}{\stepsize \varepsilon^2},
	\]
	provided that the stepsize satisfies $\stepsize \leq  \bar \stepsize/2$ for $\bar \stepsize$ such that for any $w\in \reals^d$, $v \mapsto  \bar \stepsize f(w-v) + \|v\|_2^2/2$ is convex.
	
	The number of iterations of an approximate augmented Moreau gradient descent~\eqref{eq:amp_app} with constant approximation error $\varepsilon/\sqrt 2$ to output a point that is $\stepsize\varepsilon$ close to a $2\varepsilon$ stationary point is at most
	\[
	K \leq \frac{4 (f(w^{(0)})- f^*)}{\stepsize \varepsilon^2},
	\]
	provided that $\stepsize \leq  \bar \stepsize/4$, $\reglag \leq \bar \reglag/2$ for $\bar \stepsize$ and $\bar \reglag$ such that $v \mapsto  \bar \stepsize f(w-v) + \|v\|_2^2/2$ and $v \mapsto \bar \reglag (f(w-v) -f(w))^2/2 + \|v\|_2^2/2$ are convex and $\varepsilon \leq 1/(2\reglag \lip_f)$ with $\lip_f$ the Lipschitz-continuity parameter of $f$.
\end{corollary}

\section{Detailed Chain Rules for Moreau Gradients}\label{app:chain_rule}
\subsection{Chain Rules for the Gradient of the Moreau Envelope}\label{ssec:app_simp_comp_moreau}
We introduce here some composition rules for the gradient of the Moreau envelope of simple compositions under various assumptions on the functions. 
\paragraph{Convex outer function}
The following chain rule is the basis of the reasoning presented in Section~\ref{sec:composition_moreau}. 
\chainrulea*
\begin{proof}
	The computation of the Moreau envelope amounts to solve
	\begin{align}\label{eq:decomp_moreau}
		\min_{v \in \reals^\dimvar} \stepsize \outfunc \circ \infunc(w-v) + \frac{1}{2} \|v\|_2^2 &  = \min_{v \in \reals^ d}  \sup_{\adjvar \in \dom (\stepsize \outfunc)^ *} \adjvar^ \top \infunc(w-v) - (\stepsize \outfunc)^ *(\adjvar) + \frac{1}{2}\|v\|_2^ 2.
	\end{align}
	For $\adjvar \in\reals^ \diminter$, we have that $w \rightarrow \adjvar^ \top \infunc(w)$ is $\|\adjvar\|_2\smooth_\infunc$ smooth. Therefore for $\adjvar$ such that $\smooth_\infunc \|\adjvar\|_2\leq \frac{1}{2}$, $v \rightarrow \adjvar^ \top \infunc(w-v) + \frac{1}{2}\|v\|_2^ 2$ is 1/2-strongly convex, see Lemma~\ref{lem:smooth_strg_cvx}. 	Since $\outfunc$ is $\lip_\outfunc$-Lipschitz continuous, $\dom (\stepsize \outfunc)^ * \subset \{ \adjvar \in \reals^ \diminter: \|\adjvar\|_2 \leq \stepsize \lip_\outfunc\}$. Therefore, for $\stepsize \lip_\outfunc \smooth_\infunc \leq \frac{1}{2}$, problem~\eqref{eq:decomp_moreau} is strongly convex in $v$ and concave in $\adjvar$ with $\dom (\stepsize \outfunc)^ *$ compact. Hence we can interchange min and max such that the problem reads
	\[
	\sup_{\adjvar \in \dom (\stepsize \outfunc)^*}\left\{ \min_{v\in \reals^ d} \left\{\adjvar^ \top \infunc(w-v) + \frac{1}{2}\|v\|_2^2 \right\} - (\stepsize \outfunc)^ *(\adjvar)\right\} = 	\sup_{\adjvar \in \reals^ \diminter} \env(\adjvar^ \top \infunc)(w) - (\stepsize \outfunc)^ *(\adjvar).
	\]
	Note that $\adjvar \rightarrow \env(\adjvar^\top \infunc)(w)$ is concave as an infimum of linear functions  in $\adjvar$. In addition, if $\outfunc$ is $\smooth_\outfunc$-smooth, $(\stepsize \outfunc)^ *$ is $1/(\stepsize \smooth_\outfunc)$ strongly convex. Therefore the above problem is strongly concave. For the solution $\adjvar^*$ of the above problem, the primal solution is given by 
	$
	\argmin_{v \in \reals^\dimvar}\left\{ {\adjvar^*}^\top \infunc(w - v) + \frac{1}{2}\|v\|_2^2\right\}.
	$
\end{proof}
By considering standard convex optimization algorithms to solve the subproblems defining the Moreau gradient we can estimate its computational complexity for a single composition as follows. 
\begin{lemma}\label{lem:cplxity_mb_single_comp}
Under the assumptions of Lemma~\ref{prop:chain_rule}, problem~\eqref{eq:dual_comput} can be solved up to $\varepsilon$ accuracy in at most $ O\left(\lip_\infunc^ 2\smooth_\outfunc/(\lip_\outfunc \smooth_\infunc) \ln^2 \varepsilon \right)$ calls to the gradients of $\outfunc$ and $\infunc$.
\end{lemma}
\begin{proof}
	Denoting $e(\adjvar) =\env( \adjvar^ \top \infunc)(w)$, consider a proximal gradient ascent to solve~\eqref{eq:dual_comput}, i.e., starting from $\adjvar^ {(0)}=0$, consider the iterations
	\begin{align}\label{eq:prox_grad}
		\adjvar^{(k+1)}&  = \argmax_{\nu\in \dom (\stepsize \outfunc)^*} e(\adjvar^ {(k)}) + \nabla e(\adjvar^ {(k)})^ \top (\nu -\adjvar^ {(k)}) - (\stepsize \outfunc)^ *(\nu) - \frac{1}{2\auxxstepsize} \|\nu -\adjvar^ {(k)}\|_2^ 2. 
	\end{align}
	These iterations are well-defined and converge for $\auxxstepsize$ sufficiently small provided that $e$ is smooth, i.e., differentiable with Lipschitz continuous gradients. 
	For $\adjvar \in \dom( (\stepsize \outfunc)^ *)$, we have that $v \rightarrow \adjvar^ \top\infunc(w-v) + \frac{1}{2}\|v\|_2^ 2$  is $1/2$-strongly convex, hence $\nabla e(\adjvar)$ is well defined and is given by
	\begin{align*}
		\nabla e(\adjvar) & = \infunc(w - v(\adjvar)) \\
		\mbox{for} \quad v(\adjvar)&  = \argmin_{v\in \reals^ d}  \left\{\adjvar^ \top \infunc(w-v) + \frac{1}{2}\|v\|_2^2 \right\}  =  \nabla\infunc(w- v(\adjvar)) \adjvar.
	\end{align*}
	For $\adjvar, \nu \in \dom( (\stepsize \outfunc)^ *)$, such that $\|\adjvar\| \leq \stepsize \lip_\outfunc$,   we have
	\begin{align*}
		\|v(\adjvar) - v(\nu) \|_2 & \leq \lip_\infunc\|\adjvar-\nu\|_2 + \|\adjvar\|_2 \smooth_\infunc \|v(\adjvar) - v(\nu)\|_2 \leq \lip_\infunc\|\adjvar-\nu\|_2 + \stepsize \lip_\outfunc  \smooth_\infunc \|v(\adjvar) - v(\nu)\|_2.
	\end{align*}
	Hence we have, using that $0 \leq \stepsize \leq 1/(2\lip_\outfunc \smooth_\infunc)$,
	\begin{align*}
		\|\nabla e(\adjvar) - \nabla e(\nu)\|_2 & \leq \lip_\infunc \|v(\adjvar) - v(\nu)\|_2 \leq \frac{\lip_\infunc^ 2}{1- \stepsize \lip_\outfunc\smooth_\infunc}\|\adjvar- \nu\|_2 \leq 2 \lip_\infunc^ 2\|\adjvar-\nu\|_2.
	\end{align*}
	Therefore,  the function $e$ is $2\lip_\infunc^ 2$ smooth. Since $(\stepsize \outfunc)^ *$ is $1/(\stepsize \smooth_\outfunc)$ strongly convex, a proximal gradient ascent~\eqref{eq:prox_grad} with step-size $\auxxstepsize = 1/(4\lip_\infunc^ 2)$ converges in $ \bigO(\lip_\infunc^ 2\stepsize  \smooth_\outfunc \ln \varepsilon) \leq\bigO\left( \lip_\infunc^ 2 \smooth_\outfunc/ (\lip_\outfunc \smooth_\infunc)\ln \varepsilon \right)$ iterations to a point $\varepsilon$ close to the solution ~\citep{nesterov2013gradient}.
	An accelerated gradient ascent converges in $O(\lip_\infunc\sqrt{\smooth_\outfunc/\lip_\outfunc \smooth_\infunc}\ln \varepsilon)$ iterations~\citep{nesterov2013gradient}.  Each step requires to compute, for $\dir =  \adjvar^ {(k)} + \auxxstepsize \nabla e(\adjvar^ {(k)}) =  \adjvar^ {(k)}  + \auxxstepsize \infunc(w- \nabla \env({\adjvar^{(k)} }^ \top\infunc)(w))$, 
	\begin{align*}
		\argmin_{\nu \in \reals^\diminter} (\stepsize \outfunc)^ *(\nu)  + \frac{1}{2\auxxstepsize} \|\dir- \nu\|_2^ 2
		&   = \dir -\auxxstepsize \argmin_{z \in \reals^ \diminter} \left\{\stepsize\outfunc(z) + \frac{\auxxstepsize}{2} \|z - \dir/\auxxstepsize\|_2^2\right\}  = \auxxstepsize \nabla \env((\stepsize/\auxxstepsize)\outfunc)(\dir/\auxxstepsize).
	\end{align*}
	Computing $\nabla \env((\stepsize/\auxxstepsize)\outfunc)(\dir)$ up to $\varepsilon$ accuracy requires 
	$\bigO((\stepsize/\auxxstepsize)\smooth_\outfunc \ln(\varepsilon)) \leq \bigO\left(\lip_\infunc^ 2\smooth_\outfunc/ (\lip_\outfunc\smooth_\infunc)\ln \varepsilon \right)$
	 iterations of gradient descent and 
	 $\bigO\left(\lip_\infunc\sqrt{\smooth_\outfunc /(\lip_\outfunc \smooth_\infunc)}\ln \varepsilon \right)$
	  iterations of accelerated gradient descent. 
	Computing $\nabla \env({\adjvar^{(k)} }^ \top \infunc)(w)$ for $\|\adjvar \|\leq  \stepsize\lip_\outfunc$ up to $\varepsilon$ accuracy requires $\bigO\left(\stepsize \lip_\outfunc \smooth_\infunc\ln \varepsilon \right) \leq \bigO(\ln \varepsilon)$ iterations of gradient descent or an accelerated gradient descent. 
	
	An approximate accelerated proximal gradient ascent with increasing accuracy of the inner computations of the oracles leads up to logarithmic factors to the same rates as an accelerated proximal gradient method, as demonstrated for example by \citet{lin2018catalyst} in the case of approximate accelerated proximal point methods. 
	Overall, computing the gradient of the Moreau envelope  using accelerated gradient descent in the outer and inner loops costs up to logarithmic factors $O(\lip_\infunc^ 2\smooth_\outfunc/\lip_\outfunc \smooth_\infunc)$.
\end{proof}
Finally, the following corollary illustrates the potential advantages of using Moreau gradients instead of gradients to minimize a single composition.  
\begin{corollary}
	Consider minimizing $\outfunc\circ \infunc$ for $\outfunc$, $\infunc$ satisfying the assumptions of Lemma~\ref{prop:chain_rule}. A gradient descent on the composition as in~\eqref{eq:gd_app} computes a $\varepsilon$-stationary point of the composition in at most $O((\smooth_\outfunc \lip_\infunc^ 2 + \smooth_g \lip_\outfunc)/\varepsilon^2)$ calls to the gradients of $\outfunc$ and $\infunc$. An approximate Moreau gradient descent computes a point that is $O(\lip_\infunc \smooth_\outfunc \varepsilon)$ close to a $\varepsilon$-stationary point in at most $O(\smooth_\outfunc \lip_\infunc^2  \ln^2(\varepsilon)/\varepsilon^2)$ iterations. 
\end{corollary}
\begin{proof}
	The maximum step-size for a gradient descent is a priori bounded by $1/(\smooth_\outfunc \lip_\infunc^ 2 + \smooth_\infunc \lip_\outfunc)$ using that $\smooth_{\outfunc\circ \infunc} \leq \smooth_\outfunc \lip_\infunc^ 2 + \smooth_\infunc \lip_\outfunc$ is an upper bound on the smoothness parameter of the composition. The first claim follows. The second claims follows from Lemma~\ref{lem:cplxity_mb_single_comp} and Lemma~\ref{lem:mbp_cvg}.
\end{proof}

\paragraph{Convex non-smooth composition}
Lemma~\ref{prop:chain_rule} is a particular case of the following lemma. Namely, the smoothness properties used in Prop.~\ref{prop:chain_rule} are used to characterize the range of possible stepsizes and the complexity of solving the sub-problems. Yet, in general, these smoothness assumptions are not necessary as long as the outer function is convex and its gradient takes values in the range of variables $\adjvar$ such that $v \rightarrow \adjvar^ \top \infunc(v) + \frac{1}{2} \|v\|_2^2 $ is convex.

\begin{lemma}
	Consider $\infunc:\reals^ \dimin \rightarrow \reals^\diminter$ and $\outfunc:\reals^\diminter \rightarrow\reals$.  Assume $\outfunc$ to be convex and for $\stepsize \geq 0$,
	\[
	\dom (\stepsize \outfunc)^* \subset \Lambda_\infunc := \{ \adjvar \in \reals^\diminter : v \rightarrow \adjvar^ \top \infunc(v) + \frac{1}{2} \|v\|_2^2 \ \mbox{is strongly convex}\}
	\]
	with $\dom (\stepsize \outfunc)^*$ compact, 
	then 
	\begin{align*}
		\nabla \env(\stepsize \outfunc\circ\infunc)(w)  & =  \argmin_{v \in \reals^\dimvar}\left\{ {\adjvar^*}^\top \infunc(w - v) {+} \frac{1}{2}\|v\|_2^2\right\}, \nonumber\\
		\mbox{where} \quad \adjvar^*  & \in \argmax_{\adjvar \in \reals^\diminter} - (\stepsize \outfunc)^*(\adjvar) {+}  \env(\adjvar^\top \infunc) (w).
	\end{align*}
\end{lemma}
\begin{proof}
	The computation of the Moreau envelope amounts to solve
	\begin{align}\nonumber
		\min_{v \in \reals^\dimvar} \stepsize \outfunc\circ \infunc(w-v) + \frac{1}{2} \|v\|_2^2 &  = \min_{v \in \reals^ d}  \sup_{\adjvar \in \dom (\stepsize \outfunc)^ *} \adjvar^ \top \infunc(w-v) - (\stepsize \outfunc)^ *(\adjvar) + \frac{1}{2}\|v\|_2^ 2.
	\end{align}
	By assumption, the min and sup can be swapped, which gives the result. 
\end{proof}
The definition of $\Lambda_\infunc$ can be derived for simple functions:
\begin{itemize}[nosep]
	\item for $\infunc(w) = Aw$, we have $\Lambda_\infunc = \reals^\diminter$,
	\item for $\infunc(w) = (\sigma_i(w))_{i=1}^ \diminter$ with $\sigma_i:\reals^\dimin \rightarrow\reals$ convex for all $i$, 
	we have $\Lambda_\infunc \supset \reals^\diminter_+ = \{\adjvar \in \reals^\diminter: \adjvar_i \geq 0, \ \forall i \in \{1, \ldots, \diminter\}\}$, the positive orthant,
	\item for $\infunc$ $L_\infunc$-smooth, $\Lambda_\infunc \supset \{\adjvar \in \reals^\diminter:  \|\adjvar\|_2 < 1/\smooth_\infunc\}$.
\end{itemize}
In general, we always have $0 \in \Lambda_\infunc$.

As an example, for $\outfunc$ strictly convex and continuously differentiable and $\infunc: w \rightarrow Aw$ linear the Moreau envelope of the composition can be expressed as
\[
\nabla \env(\outfunc\circ A) = A^\top\circ\left(AA^\top + \nabla \outfunc^{-1}\right)^{-1}\circ A,
\] 
which is a reformulation of~\citep[Proposition 23.25]{bauschke2011convex}.

\section{Helper Lemmas}\label{app:helper}
The following classical lemma presents that a smooth function can naturally satisfy condition~\eqref{eq:mor_cond}.
\begin{lemma}\label{lem:smooth_strg_cvx}
	If $f: \reals^ d \rightarrow \reals$ is $\smooth$-smooth, then for any $\tau\geq 1$,  $f + \frac{\smooth \tau}{2} \|\cdot\|_2^ 2$ is $(\tau-1) \smooth$ strongly convex.
\end{lemma}
\begin{proof}
	We have by smoothness of $f$, for any $w, v \in \reals^ d$, 
	$
	f(v) \geq  f(w) + \nabla f(w)^ \top (v-w) - \frac{\smooth}{2}\|v-w\|_2^ 2.
	$
	Hence, 
	\[
	f(v) + \frac{\smooth \tau}{2}\|v\|_2^2 \geq f(w) + \frac{\smooth \tau}{2} \|w\|_2^ 2 + (\nabla f(w) + \smooth \tau w)^ \top (v-w) + (\tau -1)\frac{\smooth}{2}\|v-w\|_2^ 2,
	\]
	where we used that $\|v\|_2^ 2 = \|w\|_2^ 2 + 2w^\top(v-w) + \|v-w\|_2^ 2$.  Therefore $f + \frac{\smooth \tau}{2} \|\cdot\|_2^ 2$ is $(\tau-1) \smooth$ strongly convex \citep{nesterov2013introductory}.
\end{proof}

As the Moreau gradients are computed approximately, convergence guarantees for an algorithm based on the Moreau gradients require to consider approximate oracles as done by~\citep{devolder2014first} in the convex case adapted here to assess convergence to stationary points. 
\begin{lemma}\label{lem:approx_grad_cvg}
	Let $f: \reals^d \rightarrow \reals$ be an $L$-smooth function. Consider an approximate gradient descent on $f$ with step size $0\leq \stepsize\leq1/L$, i.e.,
	$
	w^{(k+1)} = w^{(k)} - \stepsize \widehat \nabla f(w^{(k)}),
	$
	where $\|\widehat \nabla f(w^{(k)}) - \nabla f(w^{(k)})\|_2 \leq \varepsilon_k$. After $K$ iterations, this method satisfies, 
	\[
	\min_{k\in\{0, \ldots, K-1\}}\|\nabla f(w^{(k)})\|_2^2 \leq  {\frac{2(f(w^{(0)})  -\min_{w \in \reals^d} f(w))}{\stepsize K}}  +  \frac{1}{K}\sum_{k=0}^{K-1} \varepsilon_k^2.
	\]
\end{lemma}
\begin{proof}
	Denote $g^{(k)} = \widehat \nabla f(w^{(k)}) - \nabla f(w^{(k)})$ for all $k\geq 0$. By $L$-smoothness of the objective, the iterations of the approximate gradient descent satisfy, using in~$(i)$ that $L\stepsize \leq 1$,
	\begin{align*}
		f(w^{(k+1)}) & \leq f(w^{(k)}) + \nabla f(w^{(k)})^\top (w^{(k+1)} - w^{(k)}) + \frac{L}{2}\|w^{(k+1)}-w^{(k)}\|_2^2 \\
		& =  f(w^{(k)}) - \stepsize \|\nabla f(w^{(k)})\|_2^2 -\stepsize \nabla f(w^{(k)})^\top g^{(k)} + \frac{L\stepsize^2}{2}\|\nabla f(w^{(k)})+ g^{(k)}\|_2^2 \\
		&  = f(w^{(k)}) - \stepsize\left(1-\frac{L\stepsize}{2}\right)\|\nabla f(w^{(k)})\|_2^2 + \frac{L\stepsize^2}{2}\|g^{(k)}\|_2^2  + \stepsize(L\stepsize - 1)\nabla f(w^{(k)})^\top g^{(k)} \\
		& \stackrel{(i)}{\leq} f(w^{(k)}) - \stepsize\left(1-\frac{L\stepsize}{2}\right)\|\nabla f(w^{(k)})\|_2^2 + \frac{L\stepsize^2}{2}\|g^{(k)}\|_2^2  + \stepsize(1- L\stepsize)\|\nabla f(w^{(k)})\|_2 \|g^{(k)}\|_2 \\
		& \leq  f(w^{(k)}) - \stepsize\left(1-\frac{L\stepsize}{2}\right)\|\nabla f(w^{(k)})\|_2^2 + \frac{L\stepsize^2}{2}\|g^{(k)}\|_2^2  + \frac{\stepsize(1- L\stepsize)}{2}(\|\nabla f(w^{(k)})\|_2^2 + \|g^{(k)}\|_2^2) \\
		& \leq f(w^{(k)}) - \frac{\stepsize}{2} \|\nabla f(w^{(k)})\|_2^2 + \frac{\stepsize}{2} \|g^{(k)}\|_2^2.
	\end{align*}
	Summing from $k=0$ to $K-1$ and rearranging the terms, we get
	\[
	\sum_{k=0}^{K-1} \|\nabla f(w^{(k)})\|_2^2  \leq  \frac{2(f(w^{(0)}) - \min_{w \in \reals^d} f(w)) }{\stepsize} + \sum_{k=0}^{K-1}\varepsilon_k^2.
	\]
	Taking the minimum of $ \|\nabla f(w^{(k)})\|_2^2$ and dividing by $K$ we get the result. 
	
\end{proof}

\begin{fact}\label{fact:prox_grad}
	Consider for $\infunc:\reals^\dimin \rightarrow\reals^{\diminter}$, $\outfunc: \reals^{\diminter} \rightarrow \reals$, $\stepsize \geq 0$, the problem
	\[
	\max_{\adjvar \in \reals^\diminter} -(\stepsize h)^*(\adjvar) + e(\adjvar) \quad \mbox{for}\  e(\adjvar) = \env(\adjvar^ \top \infunc)(x).
	\]
	A proximal gradient step from 0 on this problem amounts to compute 
	\[
	\hat  \adjvar = \argmax_{\adjvar \in \reals^ k} \nabla e(0)^ \top \adjvar - ( \stepsize \outfunc)^ *(\adjvar) - \frac{1}{2\auxstepsize}\|\adjvar-0\|_2^ 2  = \auxstepsize \nabla \env(\auxstepsize^{-1}\stepsize \outfunc)(\infunc(x)).
	\]
\end{fact}
\begin{proof}
	We have that $e(0) = \min_{v \in \reals^\dimin} 0^\top \infunc (w - v) + \frac{1}{2}\|v\|_2^2$. So the minimizer is given as $v=0$ and $\nabla e(0) = \infunc(w)$. So the proximal gradient step amounts to solve
	\begin{align*}
		\min_{\adjvar \in \reals^\diminter}  ( \stepsize \outfunc)^ *(\adjvar) + \frac{1}{2\auxstepsize}\|\adjvar-\auxstepsize  \infunc(w)\|_2^ 2   
		& = \min_{\adjvar \in \reals^\diminter} \max_{z \in \reals^\diminter} z^\top \adjvar - \stepsize \outfunc(z)  +  \frac{1}{2\auxstepsize}\|\adjvar-\auxstepsize  \infunc(w)\|_2^ 2 \\ 
		& = \max_{z \in \reals^\diminter} - \stepsize \outfunc(z) - \frac{\auxstepsize}{2} \|z -  \infunc(w)\|_2^2.
	\end{align*}
	Hence the maximum in $z$ is reached for $\hat z = \prox{\outfunc}(\auxstepsize^{-1} \stepsize \outfunc)( \infunc(w))$ and the corresponding optimal $\adjvar$ is $\hat \adjvar = \auxstepsize  f(w) - \auxstepsize \hat z = \auxstepsize \nabla \env(\auxstepsize^{-1} \stepsize \outfunc)( f(w))$.
\end{proof}

\section{Experimental Details}\label{app:exp}
\subsection{Experimental Settings}
\paragraph{Nonlinear control of a swinging pendulum}
We consider the control of a pendulum to make it swing up after a finite time. A pendulum is described by the angle of the rod $\theta$ with the vertical axis, and its dynamics are described in continuous time as
\[
\ddot \theta(t) = -{g}\sin \theta(t)/l - {\mu} \dot\theta(t)/{ml^2} + w(t)/{ml^2},
\]
where $m=1$ denotes the mass of the bob, $l=1$ denotes the length of the rod, $\mu=0.01$ is the friction coefficient, $g=9.81$ is the $g$-force and $w(t)$ is a torque applied to the pendulum. The state of the pendulum is given by the angle $\theta$ and its speed $\omega=\dot \theta$ concatenated in $x= (\theta, \omega)$ which after discretization follow the dynamics $\phi(w_t, x_{t-1}) = x_t = ( \theta_t, \omega_t)$ s.t.
\begin{align*}
	\theta_t  & = \theta_{t-1} + \delta  \omega_{t-1} \\
	\omega_t  & = \omega_{t-1} + \delta \left(-{g}\sin \theta_{t-1}/{l} - {\mu} \omega_{t-1}/{ml^2} + w_t/{ml^2}\right),
\end{align*}
where $\delta$ is the discretization step, $x_{t-1} = (\theta_{t-1}, \omega_{t-1})$ is the current state and $w_t$ is a control parameter.
The objective is then given by $h(x_\horizon) = (\theta_\horizon- \pi)^ 2 + \rho \omega_\horizon^ 2$ for $x_\horizon= (\theta_\horizon, \omega_\horizon)$,
where $\rho=0.1$ is a penalty parameter.
The objective $h$ enforces the pendulum to swing up and be close to equilibrium (low speed) after $\horizon$ steps.

Overall the problem  can be written in the form 
\begin{align*}
	\min_{\ctrl_1, \ldots, \ctrl_\horizon} \quad &  h(\state_\horizon)  \quad \\
	\mbox{s.t.} \quad  & \state_{t+1} = \phi_t(\ctrl_t, \state_{t-1})  \ \mbox{for} \ t\in \{1, \ldots, \horizon\}, 
\end{align*}
for $\state_0$ fixed as in Eq.~\eqref{eq:chain}. 
The horizon $\horizon$ is usually large to ensure that the discretization scheme is accurate enough. In the experiments, we take a discretization step $\delta =0.1$ and a horizon $\tau = 50$ or $\tau= 100$. As many compositions are involved, we are interested in the effects of using approximate Moreau gradients compared to classical gradients.

\paragraph{Image classification with deep networks}
We consider the classifications of images from the image classification dataset CIFAR10 composed of 50 000 images classified in 10 classes~\citep{krizhevsky2009learning}. For this task, we consider either (i) a Multi-Layer Perceptron (MLP) with hidden layers $(4000, 1000, 4000)$, (ii) a convolutional network (ConvNet) with an architecture specified as
\begin{align*}
&	\mbox{  Conv[16×32×32]→ReLU→Pool[16×16×16]} \\
&	\mbox{→Conv[20×16×16]→ReLU→ Pool[20×8×8]}\\
&	\mbox{→Conv[20×8×8]→ReLU→Pool[20×4×4]}\\
&	\mbox{→FC}
\end{align*}
as done by \citet{frerix2018proximal}, where $\operatorname{Conv}[C, H, W]$ and $\operatorname{Pool}[C, H, W]$ stand for a convolutional layer and an average pooling layer respectively outputting images with $C$ channels, a height $H$  and a width $W$, (iii) a deeper convolutional network, namely, the AllCNN-C architecture presented by~\citet{springenberg2014striving}. 
On top of these architectures, we consider a cross-entropy loss and add a square regularization term $\reg(w) = \sum_{t=1}^\horizon \mu \|w_t\|_2^2/2$ with $\mu=10^{-6}$.

\subsection{Implementation details}
\paragraph{Inner computations}
The approximate back-propagation of the Moreau envelope  outlined in Algo.~\ref{algo:mbp_backward} involves an inexact minimization to access approximations of the Moreau gradients of the intermediate computations. To illustrate the potential of the proposed approach, we report the results when this inexact minimization is performed with 2 steps of a quasi-Newton algorithm. Namely, we use one step of gradient descent with a Goldstein line-search followed by one step of gradient descent using the Barzilai-Borwein step-size computation~\citep{bonnans2006numerical}.

\paragraph{Hyper-parameters}
In Algo.~\ref{algo:mbp_backward}, we choose $\scaling_t = \scaling^{\tau - t+1}$ and $\stepsize_t = \stepsize \scaling_{t+1}$ such that the updates of our algorithm can be rewritten $ \costate_{t-1} =    \widehat\nabla \env(\scaling \adjvar_{t}^ \top \phi_t(\ctrl_t, \cdot))(\state_{t-1})$ and $g_t = \widehat \nabla \env(\stepsize\adjvar_{t}^ \top \phi_t(\cdot, \state_{t-1}))(\ctrl_t)$. This choice of stepsize is motivated by Lem.~\ref{lem:stepsize} that shows that the scaling parameters $\scaling_t$ required for the subproblems to be strongly convex need to decrease geometrically as $t$ goes from $\horizon$ to $1$. 

For the nonlinear control example, we set $\scaling = 0.5$ and perform a grid search on powers of 2 for $\stepsize$ which gives $\stepsize = 2^7$ for both $\tau=50$ or $\tau=100$.  
In comparison, a grid search on the stepsize, denoted $\stepsize$, of a gradient descent on powers of 2 gives $\stepsize = 1$ for $\tau=50$ and $\stepsize=0.25$ for $\tau=100$.

For the deep learning example, we found $\scaling =1, \stepsize = 2$ after a grid-search on these parameters and $\scaling=0.5, \stepsize = 0.5$ for the ConvNet experiment. The stepsize of SGD was optimized on a grid of powers of 2 and gave $\stepsize=1$ for the MLP experiment and $\stepsize=0.125$ for the ConvNet. 

\subsection{Experimental Results}
\paragraph{Nonlinear control of a swinging pendulum}
In Fig.~\ref{fig:ctrl}, we compare a gradient descent to our Moreau gradient descent whose oracles are given by Algo.~\ref{algo:mbp_backward} on the control of a pendulum  for various horizons $\tau$.
We observe that our approach can provide smoother and faster optimization in this setting.

\paragraph{Supervised classification with deep networks}
For supervised classification with deep networks, we consider a mini-batch stochastic counterpart to the proposed algorithm as explained in Sec.~\ref{sec:oracle}. 
In Fig.~\ref{fig:deep}, we compare plain mini-batch stochastic gradient descent against a mini-batch approximate stochastic Moreau gradient descent as presented in Sec.~\ref{sec:oracle}.
The plots present the minimum of the loss or the test  error  obtained so far, i.e., on the y-axis we plot $y_k = \min_{i=0, \ldots, k} \obj(\chain(\ctrls^{(i)}))$, where $\outfunc \circ\infunc$ denote the overall test loss and $\ctrls^{(i)}$ is the current set of parameters.  
We observe that the mini-batch stochastic counterpart of the proposed algorithm compares favorably with stochastic gradient descent in both cases.

\paragraph{Stochastic algorithms with momentum}
Moreau gradients define first-order oracles that can be incorporated in popular algorithms for stochastic training with momentum such as Adam~\cite{kingma2015adam}. We illustrate this by considering the Proximal BackPropagation algorithm of~\cite{frerix2018proximal} which can be seen as a particular implementation of a Moreau gradient and apply it to the image classification dataset CIFAR10 using the AllCNN-C architecture~\cite{springenberg2014striving} with a logistic loss\footnote{A similar experiment was done by~\cite{frerix2018proximal} on a smaller architecture.}. In Fig.~\ref{fig:deep}, we observe that an approach using Moreau gradients can optimize faster on the training loss, while an approach with classical gradients can generalize better in this experiment.

\end{document}